\documentclass[12pt]{article}
\usepackage{amsmath,amssymb,amscd,array}

\let\ssection=\section
\renewcommand{\section}{\setcounter{equation}{0}\ssection}

\def\d{\delta}
\def\G{\Gamma}
\def\om{\omega}
\def\r{\rho}
\def\a{\alpha}
\def\b{\beta}
\def\s{\sigma}
\def\vfi{\varphi}

\def\l{\lambda}
\def\m{\mu}
\def\n{\nabla}

\def\implies{\Rightarrow}

\newcommand{\bbR}{\mathbb{R}}

\newcommand{\Og}{\mathrm{O}}
\newcommand{\Diff}{\mathrm{Diff}}

\newcommand{\Div}{\mathrm{Div}}
\newcommand{\End}{\mathrm{End}}
\newcommand{\cF}{{\mathcal{F}}}
\newcommand{\cD}{{\mathcal{D}}}

\newcommand{\gl}{{\mathrm{gl}}}

\newcommand{\PSL}{\mathrm{PSL}}
\newcommand{\SL}{\mathrm{SL}}
\newcommand{\Sl}{\mathrm{sl}}

\newcommand{\Vect}{\mathrm{Vect}}

\newcommand{\og}{\mathrm{g}}

\newcommand{\cqfd}{\hspace*{\fill}\rule{3mm}{3mm}}
\newcommand{\cqf}{\hspace*{\fill}\rule{2mm}{2mm}}

\begin{document}

\frenchspacing

\def\d{\delta}
\def\g{\gamma}
\def\om{\omega}
\def\r{\rho}
\def\a{\alpha}
\def\b{\beta}
\def\s{\sigma}
\def\vfi{\varphi}
\def\l{\lambda}
\def\m{\mu}
\def\implies{\Rightarrow}

\oddsidemargin .1truein
\newtheorem{thm}{Theorem}[section]
\newtheorem{lem}[thm]{Lemma}
\newtheorem{cor}[thm]{Corollary}
\newtheorem{pro}[thm]{Proposition}
\newtheorem{ex}[thm]{Example}
\newtheorem{rmk}[thm]{Remark}
\newtheorem{defi}[thm]{Definition}
\title{Projective and Conformal Schwarzian Derivatives and
Cohomology of Lie Algebras Vector Fields Related to Differential
Operators}
\author{Sofiane Bouarroudj\\
{\footnotesize Department of Mathematics, U.A.E. University,
Faculty of Science} \\
{\footnotesize P.O. Box 15551, Al-Ain, United Arab Emirates.}\\
{\footnotesize  e-mail:bouarroudj.sofiane@uaeu.ac.ae} }
\date{}
\maketitle
\newpage
\begin{abstract} Let $M$ be either a projective manifold $(M,\Pi)$
or a pseudo-Riemannian manifold $(M,\og).$ We extend,
intrinsically, the projective/conformal Schwarzian derivatives
that we have introduced recently, to the space of differential
operators acting on symmetric contravariant tensor fields of any
degree on $M.$ As operators, we show that the projective/conformal
Schwarzian derivatives depend only on the projective connection
$\Pi$ and the conformal class $[\og]$ of the metric, respectively.
Furthermore, we compute the first cohomology group of $\Vect(M)$
with coefficients into the space of symmetric contravariant tensor
fields valued into $\delta$-densities as well as the corresponding
relative cohomology group with respect to $\Sl(n+1,\mathbb{R}).$
\end{abstract}
\section{Introduction}
The investigation of invariant differential operators is a famous
subject that have been intensively investigated by many authors.
The well-known invariant operators and more studied in the
literature are the {\it Schwarzian derivative}, the power of the
Laplacian (see \cite{eas}) and the Beltrami operator (see
\cite{be}). We have been interested in studying the Schwarzian
derivative and its relation to the geometry of the space of
differential operators viewed as a module over the group of
diffeomorphisms in the series of papers \cite{bp1,bo2,bo1}. As a
reminder, the classical expression of the Schwarzian derivative of
a diffeomorphism $f$ is:
\begin{equation}
\frac{f'''}{f'}-\frac{3}{2}\left (\frac{f''}{f'}\right )^2\cdot
\label{lam}
\end{equation}
The two following properties of the operator (\ref{lam}) are the
most of interest for us:

(i) It vanishes on the M\"{o}bius group $\PSL(2,\bbR)$ -- here the
group $\SL(2,\bbR)$ acts locally on $\mathbb{R}$ by projective
transformations.

(ii) For all diffeomorphisms $f$ and $g,$ the equality
\begin{equation}
\label{lam2}
 S(f\circ g)={g'}^2\cdot S(f)\circ g+S(g)
\end{equation}
holds true.

The equality (\ref{lam2}) seems to be known since Cayley; however,
it was first reported by Kirillov and Segal (see \cite{ki,ko,se})
that this property is nothing but a 1-cocycle property -- it
should be stressed that cocycles on the group are not easy to come
up with, and only few explicit expressions are known (cf.
\cite{f}).

Our study has its genesis from the geometry of the space of
differential operators acting on tensor densities, viewed as a
module over the group of diffeomorphisms and also over the Lie
algebra of smooth vector fields. In the one-dimensional case, this
study have led to compute the (relative) cohomology group
$$
\mathrm{H^1}(\Diff(\mathbb{R}),\PSL(2,\bbR);
\End_{\mathrm{diff}}(\cF_\lambda,\cF_\mu)),
$$
where $\cF_\lambda$ is the space of tensor densities of degree
$\lambda$ on $\mathbb{R}.$

 It turns out that the Schwarzian derivative as well as
new cocycles span the cohomology group above, as proved in
\cite{bo1}. These new 1-cocycles can also be considered as natural
generalizations of the Schwarzian derivative (\ref{lam}), although
they are only defined on an one-dimensional manifold.

The first step towards generalizing the Schwarzian derivative
underlying the properties (i) and (ii) to multi-dimensional
manifolds was a part of our thesis \cite{hm}. It was aimed at
defining the {\it projective} Schwarzian derivatives as 1-cocycles
on $\Diff(\mathbb{R}^n)$ valued into the space of differential
operators acting on contravariant twice-tensor fields, and vanish
on $\PSL(n+1,\bbR).$ Later on, we constructed in \cite{bp1}
1-cocycles on $\Diff(\mathbb{R}^n)$ valued into the same space but
vanish on the conformal group $\Og(p+1,q+1),$ where $p+q=n.$ These
$\Og(p+1,q+1)$-invariant 1-cocycles were interpreted as {\it
conformal} Schwarzian derivatives. Moreover, these
projectively/conformally invariant 1-cocycles were built
intrinsically by means of a projective connection and a
pseudo-Riemannian metric, thereby making sense on any curved
manifold. As projective structures and conformal structures
coincide in the one-dimensional case, these (projective/conformal)
1-cocycles are considered as natural generalizations of the
Schwarzian derivative (\ref{lam}).

This paper is, first, devoted to extend these derivatives to the
space of differential operators acting on symmetric contravariant
tensor fields of any degree.

In virtue of the one-dimensional case, the (projective/conformal)
Schwarzian derivatives should define cohomology classes belonging
to
$$
\mathrm{H^1}(\Diff(\mathbb{R}^n),\mathfrak{H};
\End_{\mathrm{diff}}\,({\cal S}_\delta(\mathbb{R}^n),{\cal
S}_\delta(\mathbb{R}^n))),
$$
where ${\cal S}_\delta(\mathbb{R}^n)$ is the space of symmetric
contravariant tensor fields on $\mathbb{R}^n$ valued into
$\delta$-densities and $\mathfrak{H}$ is the Lie group
$\PSL(n+1,\mathbb{R})$ or $\Og(p+1,q+1).$

The cohomology group above is not easy to handle; nevertheless, we
compute in Theorem \ref{sala} the cohomology group
$$
\mathrm{H^1}(\Diff(\mathbb{S}^n); \End_{\mathrm{diff}}\,({\cal
S}_\delta(\mathbb{S}^n),{\cal S}_\delta(\mathbb{S}^n))),
$$
for the (two and three)-dimensional sphere.

Moreover, we compute in Theorem \ref{class} the (relative)
cohomology group
\begin{equation}
\label{nco} \mathrm{H^1}(\Vect(\mathbb{R}^n),\Sl(n+1,\bbR);
\End_{\mathrm{diff}}\,({\cal S}_\delta(\mathbb{R}^n),{\cal
S}_\delta(\mathbb{R}^n)))\cdot
\end{equation}
The computation being inspired from Lecomte-Ovsienko's work
\cite{lo2}, uses the well-known Weyl's classical invariant theory
\cite{Weyl}. It provides a proof -- at least in the infinitesimal
level -- that the {\it infinitesimal} projective Schwarzian
derivatives that we are introducing are unique.

Furthermore, we compute in Theorem \ref{tah} the cohomology group
$$
\mathrm{H^1}(\Vect(M); \End_{\mathrm{diff}}\,({\cal
S}_\delta(M),{\cal S}_\delta(M))),
$$
where $M$ is an arbitrary manifold.

According to the Neijenhuis-Richardson's theory of deformation
\cite{nr}, the cohomology group above will measure all {\it
infinitesimal} deformations of the $\Vect(M)$-module ${\cal
S}_\delta (M).$
\section{The space of symbols as modules over
$\Diff(M)$ and $\Vect(M)$}
Throughout this paper, $M$ is an (oriented) manifold of dimension
$n$ endowed with an affine symmetric connection. We denote by
$\Gamma$ the Christoffel symbols of this connection and by
$\nabla$ the corresponding covariant derivative. It should be
clear from the context wether the connection is arbitrary or a
Levi-Civita one.

We use the Einstein convention summation over repeated indices.

Our symmetrization does not contain any normalization factor.
\subsection{The space of tensor densities}
The space of tensor densities of degree $\d$ on $M$, denoted by
$\cF_{\d}(M),$ is the space of sections of the line bundle:
$|\wedge^n T^*M |^{\otimes \d},$ where $\delta\in \mathbb{R}.$ In
local coordinates $(x^i),$ any $\delta$-density can be written as
$$
\phi (x)\,|dx^1\wedge\cdots \wedge dx^n|^{\delta}\cdot
$$
As examples, $\cF_{0}(M)=C^{\infty}(M)$ and
$\cF_{1}(M)=\Omega^1(M).$

The affine connection $\Gamma$ can be naturally extended to a
connection that acts on $\cF_{\d}(M).$ The covariant derivative of
a density $\phi \in \cF_{\d}(M)$ is given as follows. In local
coordinates $(x^i),$ we have
$$
\nabla_i\,\phi=\partial_i\phi-\delta\, \Gamma_{iu}^u\, \phi,
$$
where $\partial_i$ stands for the partial derivative with respect
to $x^i.$
\subsection{The space of tensor fields as a module}
Denote by ${\cal S} (M)$ the space of contravariant symmetric
tensor fields on $M.$ This space is naturally a module over the
group $\Diff(M)$ by the natural action. Moreover, it is isomorphic
to the space of symbols, namely functions on the cotangent bundle
$T^*M$ that are polynomial on fibers.

We are interested in defining a one-parameter family of
$\Diff(M)$-modules on ${\cal S} (M)$ by
$$
{\cal S}_{\d}(M):={\cal S}(M)\otimes \cF_{\d}(M).
$$
The action is defined as follows. Let $f\in \Diff(M)$ and $P\in
{\cal S}_{\d}(M)$ be given. Then, in a local coordinates $(x^i)$,
we have
\begin{eqnarray}
\label{actsym}
f_{\d}^* P&=& f^*P\cdot (J_{f^{-1}})^{\d},
\end{eqnarray}
where $J_f=|Df/Dx|$ stands for the Jacobian of $f,$ and $f^*$
stands for the natural action of $\Diff(M)$ on ${\cal S}(M).$

By differentiating the action (\ref{actsym}) we get the
infinitesimal action of $\Vect(M):$ for all $X\in \Vect(M),$ and
for all $P \in {\cal S}(M)$ we have
\begin{equation}
\label{Lieactsym} L_X P=L_X(P)+\delta\, \Div X\,P,
\end{equation}
where $\Div$ is the divergence operator associated with some
orientation.

Denote by ${\cal S}^k_\d (M)$ the space of symmetric tensor fields
of degree $k$ on $M$ endowed with the $\Diff(M)$-module structure
(\ref{actsym}). We then have a graduation of $\Diff(M)$-modules:
${\cal S}_{\d}(M)=\oplus_{k\geq 0}{\cal S}^k_{\d}(M).$

The actions (\ref{actsym}) and (\ref{Lieactsym}) are of most
interest of us. Throughout this paper, all actions will be refered
to them.
\section{A compendium on projective and conformal structures}
\label{compo}
We will collect, in this section, some gathers on projective and
conformal structures. These notions are well-known in projective
and conformal geometry. However, they are necessary to introduce
here in order to write down explicit expressions of the Schwarzian
derivatives.
\subsection{Projective structures}
A {\it projective connection} is an equivalent class of symmetric
affine connections giving the same non-parameterized geodesics.

Following \cite{kn}, the symbol of the projective connection is
given by the expression
\begin{equation}
\label{connection}
\Pi_{ij}^k=\G_{ij}^k-\frac{1}{n+1}\left
(\delta_i^k \G_{lj}^l+\delta_j^k \G_{il}^l\right ).
\end{equation}
Two affine connections $\Gamma$ and $\tilde \Gamma$ are {\it projectively
equivalent} if the corresponding symbols (\ref{connection}) coincide.
Equivalently, if there exists a 1-form $\omega$ such that
\begin{equation}
\label{assoc} \tilde
\Gamma_{ij}^k=\Gamma_{ij}^k+\delta_{j}^k\,\omega_i+\delta_{i}^k\,\omega_j.
\end{equation}
A projective connection on $M$ is called \textit{flat} if in a
neighborhood of each point there exists a local coordinates such
that the symbols $\Pi_{ij}^k$ are identically zero (see \cite{kn}
for a geometric definition).

A projective structure on $M$ is given by a local action of the
group $\SL(n+1,\bbR)$ on it. Every flat projective connection
defines a projective structure on~$M$.

On $\bbR^n$ with its standard projective structure, the Lie
algebra $\Sl(n+1,\bbR)$ can be embedded into the Lie algebra
$\Vect(\bbR^n)$ by
\begin{equation}
\label{pr} \frac{\partial}{\partial x^i},\quad
x^i\frac{\partial}{\partial x^j},\quad
x^ix^k\frac{\partial}{\partial x^k},\quad i,j=1,\ldots,n.
\end{equation}
where $(x^i)$ are the coordinates of the projective structure. The
first two vector fields form a Lie algebra isomorphic to the
affine Lie algebra $\gl(n,\mathbb{R})\ltimes \mathbb{R}^n.$
\subsection{Conformal structures}

A {\it conformal} structure on a manifold is an equivalence class
of pseudo-Riemannian metrics $[\og]$ that have the same direction.

If $\Gamma_{ij}^k$ are the Levi-Civita connection associated with
the metric $\og$, then the Levi-Civita connection, $\tilde
\Gamma_{ij},$ associated with the metric $e^{2F}\cdot \og,$ where
$F$ is a function on $M$, are related, in any local coordinates
$(x^i),$ by
\begin{equation}
\label{tof}
\tilde \Gamma_{ij}^k=\Gamma^k_{ij} +F_i
\,\delta^k_j+F_j \,\delta^k_i - \og_{ij}\,\og^{kt}\,F_t,
\end{equation}
where $F_i=\partial F/\partial x^i.$

A conformal structure on $(M,\og)$ is called \textit{flat} if in a
neighborhood of each point there exists a local coordinate system
such that the metric $\og$ is a multiple of $\og_0,$ where $\og_0$
is the metric $\mathrm{diag}(1,\ldots,1,-1,\ldots,-1)$ whose trace
is $p-q.$

It is well-know that the group of diffeomorphisms of $\bbR^n$ that
keep the standard metric $\og_0$ in the conformal class is the
group $\Og(p+1,q+1),$ where $p+q=n.$
\begin{rmk} {\rm
The Lie algebra $\mathrm{o}(p+1,q+1)$ can also be embedded into
$\Vect(\mathbb{R}^n)$ via formulas analogous to (\ref{pr}) but we
do not need them here.}
\end{rmk}
\subsection{An intrinsic 1-cocycle and a Lie derivative of
a connection}
A connection itself is not a well-defined geometrical object.
However, the difference between two connections is a well-defined
tensor fields of type $(2,1).$ Therefore, the following object
\begin{equation}
\label{tenso} \mathfrak{L}(f):=f^*\Gamma-\Gamma,
\end{equation}
where $f$ is a diffeomorphism, is globally defined on $M$.

It is easy to see that the map
$$
f\mapsto \mathfrak{L}(f^{-1})
$$ defines a
1-cocycle on $\Diff(M)$ with values into tensor fields of type
$(2,1).$

The infinitesimal 1-cocycle associated with the tensor
(\ref{tenso}), denoted by $\mathfrak l,$ is called the {\it Lie
derivative} of a connection; it can also be defined as follows.
For all $X\in \Vect(M),$ the 1-cocycle ${\mathfrak l}(X)$ is the
map
\begin{equation}
\label{mou} (Y,Z)\mapsto [X,\n_YZ]-\n_{[X,Y]} Z- \n_{Y} [X,Z]\cdot
\end{equation}

We will use intensively, throughout this paper, the tensor
(\ref{tenso}) as well as the tensor (\ref{mou}).
\section{Projectively invariant Schwarzian derivatives}
Let $\Pi$ and $\tilde \Pi$ be two projective connections on~$M.$
Then the difference $\Pi-\tilde \Pi$ is a well-defined
$(2,1)$-tensor field. Therefore, it is clear that a projective
connection on $M$ leads to the following 1-cocycle on $\Diff(M)$:
\begin{equation}
{\mathfrak T}(f^{-1})^k_{ij}:= (f^{-1})^*\Pi_{ij}^k-\Pi_{ij}^k,
\label{ell}
\end{equation}
which vanishes on (locally) projective diffeomorphisms.
\begin{rmk}{\rm There is also an alternative approach in defining
the 1-cocycle (\ref{ell}) by means of the tensor (\ref{tenso}). }
\end{rmk}
\subsection{The main definitions}
\begin{defi}{\rm
For all $f\in\Diff(M)$ and for all $P\in {\cal S}_\delta^k(M),$ we
put
\begin{equation}
\label{Scw2}
 {\mathfrak U}(f)\,(P)^{i_1\cdots i_{k-1}}=\sum_{s=1}^{k-1}{\mathfrak
T}(f)^{i_s}_{ij}\,\,P^{iji_1\cdots \widehat{i}_s\cdots i_{k-1}},
\end{equation}
where ${\mathfrak T}(f)$ is the tensor (\ref{ell}).}
\end{defi}

By construction, the operator (\ref{Scw2}) is projectively
invariant, viz it depends only on the projective class of the
connection.
\begin{thm}
\label{thm1} (i) For all $\delta\not=\frac{2k-1+n}{1+n},$ the map
$f\mapsto {\mathfrak U}(f^{-1})$ defines a non-trivial 1-cocycle
valued into
${\cal D}({\cal S}^{k}_\delta(M),{\cal S}^{k-1}_\delta(M));$\\
(ii) for $\delta=\frac{2k-1+n}{1+n},$ we have
$$
\mathfrak{U}(f)(P)^{i_1\cdots i_{k-1}}=\left ({f^{-1}}^*\,
\nabla_j-\nabla_j\right )P^{ji_1\cdots i_{k-1}}\cdot
$$
\end{thm}
{\bf Proof.} (i) The 1-cocycle property of the operator
(\ref{Scw2}) follows immediately from the 1-cocycle property of
the tensor (\ref{ell}). Let us prove the non-triviality. Suppose
that there exists an operator $A$ such that
\begin{equation}
\label{tri}
{\mathfrak U}(f)={f^{-1}}^*A-A.
\end{equation}
As $\mathfrak{U}(f)$ is a zero-order operator, the operator $A$ is
almost first-order. If $A$ is zero-order, namely a multiplication
operator, its principal symbol, say $a,$ transforms under
coordinates change as a tensor fields of type $(2,1).$ The
equality above implies that ${\mathfrak T}(f)={f^{-1}}^*a-a$ which
is absurd, as $\mathfrak T$ is a non-trivial 1-cocycle. Suppose
then that $A$ is a first-order operator, namely
$$
A(P)^{i_1\cdots i_{k-1}}=\nabla_j\,P^{ji_1\cdots i_{k-1}},
$$
for all $P\in {\cal S}^{k}_\delta(M).$  It is a matter of direct
computation to prove that
$$
({f^{-1}}^*A-A)(P)^{i_1\cdots
i_{k-1}}=\sum_{s=1}^{k-1}\mathfrak{L}(f)_{ij}^{i_s}\,
P^{iji_1\cdots \widehat{i}_s\cdots
i_{k-1}}-(\delta-1)\,\mathfrak{L}(f)_j\, P^{ji_1\cdots i_{k-1}},
$$
where $\mathfrak{L}(f)^k_{ij}$ are the components of the tensor
(\ref{tenso}). We can easily seen that the equality (\ref{tri})
holds true if and only if $\delta=\frac{2k-1+n}{1+n}.$

We will introduce a second 1-cocycle valued into ${\cal D}({\cal
S}^k(M),{\cal S}^{k-2}(M)).$ But, at first, we start by giving its
expression when $k=2.$
\begin{defi}{\rm
For all $f\in\Diff(M)$ and for all $P\in {\cal S}_\delta^2(M),$ we
put
\begin{gather}
\mathfrak{V}(f)(P):=\mathfrak{T}(f)^k_{ij} \n_k P^{ij}+
\n_k\,\mathfrak{L}(f)^k_{ij}\,\,P^{ij}-
\frac{3+n-\delta(1+n)}{1+n}\,\n_i\, \mathfrak{L}(f)_{j}\,P^{ij}
\label{MultiSchwar2}\\
\phantom{\mathfrak{V}(f):=}{}+(1-\delta)\left
(\mathfrak{L}(f)^u_{ij}\mathfrak{L}(f)_u -\frac{1}{n+1}\,
\mathfrak{L}(f)_{i}\mathfrak{L}(f)_{j}+\frac{1+n}{n-1}\left (
{f^*}^{-1}R_{ij}-R_{ij}\right)\right )P^{ij},\nonumber
\end{gather}
where $\mathfrak{L}(f)_{ij}^k$ are the components of the 1-cocycle
(\ref{tenso}), $\mathfrak{T}_{ij}^k(f)$ are the components of the
1-cocycle (\ref{ell}) and $R_{ij}$ are the components of the Ricci
tensor. }
\end{defi}
\begin{thm}
\label{thm2} (i) For all $\delta \not= \frac{n+2}{n+1}$, the map
$f\mapsto {\mathfrak V} (f^{-1})$ defines a non-trivial 1-cocycle
on $\Diff(M)$ with values into ${\cal D} ({\cal S}^2_{\d}(M),{\cal
S}^0_\d(M))$.

(ii) For $\delta = \frac{n+2}{n+1}$, we have
$$ {\mathfrak V}(f)={f^{-1}}^*B-B,$$
where $B$ is the operator
\begin{equation}
\label{yam} B:=\nabla_i\nabla_j-\frac{1}{n-1}\,R_{ij}.
\end{equation}
 (ii) The operator (\ref{MultiSchwar2}) depends only on the
projective class of the connection. When $M=\bbR^n$ (or $M=S^n$)
and $M$ is endowed with a flat projective structure, this operator
vanishes on the projective group $\PSL(n+1,\bbR).$
\end{thm}
\begin{rmk}{\rm (i)
The operator ${\mathfrak V}(f)$ in (\ref{MultiSchwar2}) enjoys the
elegant expression:
\begin{equation}
\label{iya}
{\mathfrak T}(f)^k_{ij} \n_k- \frac{2-\d
(n+1)}{n-1}\,\n_k\left({\mathfrak T} (f)^k_{ij}\right)+
\frac{(n+1)(1-\d)}{n-1}\,{\mathfrak T} (f)^k_{im}{\mathfrak
T}(f)^m_{kj},
\end{equation}
which can be obtained through the relation
\begin{equation}
\label{ric} {f^{-1}}^* R_{jk}-R_{jk}= -\n_i\,
\mathfrak{L}(f)^i_{jk}+\n_j \, \mathfrak{L}(f)_k +
\mathfrak{L}(f)^m_{sj}\, \mathfrak{L}(f)^s_{km}- \mathfrak{L}(f)_m
\, \mathfrak{L}(f)^m_{jk}.
\end{equation}
(ii) We will retain the Ricci tensor into the explicit expression
of the Schwarzian derivatives disregarding the equation
(\ref{ric}), because it will be useful when we will study theirs
relation to the well-known Vey cocycle.}
\end{rmk}
For $k>2,$ we state the following definition.
\begin{defi}
{\rm For all $f\in \Diff(M),$ and for all $P\in {\cal
S}^k_\delta(M),$ we put
\begin{gather}
\nonumber {\mathfrak V}(f)\,(P)^{i_1\cdots i_{k-2}}=
\sum_{s=1}^{k-2}\,{\mathfrak T}(f)^{i_s}_{tu}\,
\,\n_{v}\,P^{tuvi_1\cdots \widehat{i}_s\cdots i_{k-2} }
+\alpha_1\,\,{\mathfrak
T}(f)^t_{uv}\, \,\n_t\,P^{uvi_1\cdots i_{k-2} }\\
\nonumber \phantom{\mathfrak{V}(f)}{} +\sum_{s=1}^{k-2}\left
(\alpha_2\n_{t}\,
\mathfrak{L}(f)_{uv}^{i_s}+\alpha_3\,\mathfrak{L}(f)_{wt}^{i_s}
\,\mathfrak{L}(f)_{uv}^w +\alpha_4\,\mathfrak{L}(f)^{i_s}_{tu}\,
\mathfrak{L}(f)_{v}\right )
P^{tuvi_1\cdots \widehat{i}_s\cdots i_{k-2} }\\
\phantom{\mathfrak{V}(f):=}{}
\nonumber+\left (\alpha_5\,\n_t
\mathfrak{L}(f)^t_{uv}+\alpha_6 \,\nabla_u\,
\mathfrak{L}(f)_v+\alpha_7\,\mathfrak{L}(f)_{u}\,\mathfrak{L}(f)_{v}+
\alpha_8\,\mathfrak{L}(f)_{uv}^w\,\mathfrak{L}(f)_{w}\right )
\,P^{uvi_1\cdots i_{k-2}}\\
\phantom{\mathfrak{V}(f):=}{} \displaystyle +\alpha_9\,\left (
{f^*}^{-1}R_{uv}-R_{uv}\right)\,P^{uvi_1\cdots
i_{k-2}}+e\!\!\!\sum_{\substack{1\le s<r\le k-2}
}^{k-2}\,\mathfrak{L}(f)_{uv}^{i_s}\,\mathfrak{L}(f)_{pq}^{i_r}\,
P^{uvpqi_1\cdots \widehat{i}_s\cdots \widehat{i}_r\ldots i_{k-2}},
\label{MultiSchwar2k}
\end{gather}
where $R_{uv}$ are the Ricci tensor components,
$\mathfrak{L}(f)^k_{ij}$ are the components of the tensor
(\ref{tenso}) and  ${\mathfrak T}(f)^k_{ij}$ are the components of
the tensor (\ref{ell}). The
constant $ e=~\begin{cases} 1& \text{if $k\geq 4$},\\
0& \text{otherwise}
\end{cases}$ and the constants $\alpha_1,\ldots,\alpha_9$ are given
by
\begin{flalign}
\nonumber \alpha_1&=\displaystyle \frac{1}{2}\,(3-2k+
n(\delta-1)+\delta);& \displaystyle \alpha_5 &= \displaystyle
\frac{1}{2}\,(3-2k+
n(\delta-1)+\delta); \\[3mm]
\nonumber \alpha_2&=\displaystyle \frac{1}{6}\,(2k +(1-\delta)\,
(1 + n)); & \displaystyle  \alpha_6 &=\displaystyle
\frac{1}{2}\,(\delta
-1)(1-2k  +n(\delta-1)+\delta);\\[3mm]
\label{const} \displaystyle \alpha_3&=\displaystyle
\frac{1}{3}\,(5-2k+ n(\delta-1)+\delta) ; &\alpha_7&=\displaystyle
\frac{1}{2}\,(
\delta-1)^2;\\[3mm]
\nonumber \alpha_4&=(1-\delta);& \alpha_8 &=\displaystyle
\frac{1}{2}\,(1-\delta)(3-2k +n(\delta-1)+\delta);&
\end{flalign}
$$
\alpha_9= \displaystyle \frac{11+4k^2+(\delta
-1)\,(2n\,(5-4k+3\delta)+3n^2(\delta-1))+10 \delta + 3\delta^2
-4k(3+2\delta)}{6-6n}\cdot$$ }
\end{defi}
\begin{thm}
\label{thm3} (i) For all $\delta \not= \frac{2k-2+n}{n+1}$, the
map $f\mapsto {\mathfrak V} (f^{-1})$ defines a non-trivial
1-cocycle on $\Diff(M)$ with values into ${\cal D} ({\cal
S}^k_{\d}(M),{\cal S}^{k-2}_\d(M))$.

(ii) The operator (\ref{MultiSchwar2k}) depends only on the
projective class of the connection. When $M=\bbR^n$ (or
$M=\mathbb{S}^n$) and $M$ is endowed with a flat projective
structure, this operator vanishes on the projective group
$\PSL(n+1,\bbR).$
\end{thm}
We will prove Theorem (\ref{thm2}) and Theorem (\ref{thm3})
simultaneously.\\
\noindent{\bf Proof Theorem (\ref{thm2}) and Theorem
(\ref{thm3}).}
To prove that the map $f\mapsto {\mathfrak V}
(f^{-1})$ is a 1-cocycle we have to verify the 1-cocycle condition
\begin{equation}
\label{cond} {\mathfrak V}(f\circ g)={g^{-1}}^*{\mathfrak
V}(f)+{\mathfrak V}(g) \quad \mbox{for all } f,g\in \Diff(M),
\end{equation}
where $g^*$ is the natural action on ${\cal D} ({\cal
S}^k_{\d}(M),{\cal S}_\d^{k-2}(M)).$ In order to prove this
condition we will, first, remove the Ricci tensor from the
expressions (\ref{MultiSchwar2}) and (\ref{MultiSchwar2k}),
because it is obviously a coboundary; secondly, we use the
equalities
\begin{eqnarray}
\label{comp} \n_u\, f^*_{\d} P^{i_1\ldots i_k}&=&f^*_{\d}\n_u
P^{i_1\ldots i_k}-\sum_{s=1}^k \left(
\mathfrak{L}(f^{-1})_{uv}^{i_s}\, f^*_{\d} \, P^{vi_1\ldots
\widehat{i}_s\ldots i_k}\right)+\d
\,\mathfrak{L}(f^{-1})_u\, f^*_{\d} P^{i_1\ldots i_k},\nonumber\\
\mathfrak{L}(f\circ
g)^u_{ij}&=&{g^*}^{-1}\mathfrak{L}(f)^u_{ij}+\mathfrak{L}(g)^u_{ij},
\end{eqnarray}
and the equality
\begin{eqnarray}
\label{sat} \nonumber \n_u g^* \mathfrak{L}(f)_{ij}^k&=&g^*\n_u
\mathfrak{L}(f)_{ij}^k - h^*\mathfrak{L}(f)_{ij}^t\,
\mathfrak{L}(g^{-1})_{ut}^k + \mathrm{Sym}_{i,j}\left(
g^*\mathfrak{L}(f)_{it}^k \, \mathfrak{L}(g^{-1})_{ju}^t \right),
\end{eqnarray}
where $\mathfrak{L}(f)_{ij}^k$ are the components of the tensor
(\ref{tenso}). The 1-cocycle condition for the operator
(\ref{MultiSchwar2k}) can verified by a long and tedious
computation. We will give a proof here only when $k=2.$ By using
the equalities above we see that, for all $P\in {\cal
S}^2_\delta(M),$ we have
$$
\begin{array}{ccl}
\nonumber {\mathfrak V}(f\circ g)(P)&= &\left
({g^*}^{-1}\,{\mathfrak T}(f)_{ij}^k+{\mathfrak T}
(g)_{ij}^k\right )\n_k  P^{ij} +\n_k \left
({g^*}^{-1}\mathfrak{L}(f)_{ij}^k+\mathfrak{L}(g)_{ij}^k
\right )P^{ij}\\[3mm]
& &\displaystyle -\frac{3+n-\delta(1+n)}{1+n}\,\n_i\left(
{g^*}^{-1}\mathfrak{L}(f)_j+\mathfrak{L}(g)_{j}\right
)P^{ij}\\[3mm]
& &\displaystyle+(1-\d)\,\left ( \mathfrak{L}(f\circ
g)_{ij}^u\,\,\mathfrak{L}(f\circ g)_u-\frac{1}{n+1} \mathfrak{L}
(f\circ g)_{i}\,\,\mathfrak{L} (f\circ g)_{j}
\right )P^{ij}\\[3mm]
&=&g^*_{\d} \left ({\mathfrak V}(f) \,{g^*_{\d}}^{-1} (P) \right
)+{\mathfrak V} (g)(P)
\end{array}
$$
Now we prove that the 1-cocycles (\ref{MultiSchwar2}) and
(\ref{MultiSchwar2k}) are not trivial. Suppose that there exists
an operator $A: {\cal S}_{\d}^k(M)\rightarrow {\cal
S}_{\d}^{k-2}(M)$ such that
\begin{equation}
\label{assile} {\mathfrak V}(f)={f^*}^{-1} A -A.
\end{equation}
Since the operators (\ref{MultiSchwar2}) and (\ref{MultiSchwar2k})
are first-order, the operator $A$ is at most second-order. If the
operator $A$ is first-order, its principal symbol should
transforms under coordinates change as a tensor fields of type
$(2,1).$ From the equality (\ref{assile}) one can easily seen that
${\mathfrak T} (f)_{ij}^k$ is a trivial 1-cocycle, which is
absurd. If $A$ is second-order, its principal symbol should be
equal to the identity, otherwise the equality (\ref{assile}) does
not hold true. Therefore, the operator $A$ is given by
$$
A(P)^{i_1\cdots i_{k-2}}=\nabla_u\,\nabla_v P^{uvi_1\cdots
i_{k-2}},
$$
for all $P\in {\cal S}_\delta^k(M).$ Now, an easy computation
gives
\begin{gather}
{f^{-1}}^*A-A= \displaystyle
\sum_{s=1}^{k-2}\mathfrak{L}(f)_{tu}^{i_s}\, \n_v\,
P^{tuvi_1\cdots \widehat{i}_s\cdots i_{k-2}}+
(1-\delta)\,\mathfrak{L}(f)_{v}\, \n_u\,
P^{uvi_1\cdots i_{k-2}}   \nonumber \\
\phantom{{f^{-1}}^*A-A}{}\displaystyle -\sum_{s=1}^{k-2}\left (
{f^*}^{-1}\,\nabla_u\,\mathfrak{L}(f^{-1})_{tv}^{i_s}\,\,
P^{tuvi_1\cdots \widehat{i}_s\cdots i_{k-2}}+
\,\mathfrak{L}(f)_{tu}^{i_s}\, {f^*}^{-1}\,\n_v\,f^*
P^{uvti_1\cdots i_{k-2}}\right)
\nonumber \\
\phantom{f}{} -\left({f^*}^{-1}\,\nabla_r\,
\mathfrak{L}(f^{-1})_{uv}^{r}\,
-(\delta-1)\,{f^*}^{-1}\,\nabla_u\,\mathfrak{L}(f^{-1})_{v}\right)
P^{uvi_1\cdots \widehat{i}_s\cdots i_{k-2}}   \nonumber \\
\phantom{{f^{-1}}^*A}{}\displaystyle
+\,\mathfrak{L}(f)_{vu}^{r}\,{f^*}^{-1}\n_r\,f^* P^{uvi_1\cdots
i_{k-2}} -(\delta-1) \,\mathfrak{L}(f)_{u}{f^*}^{-1}\n_v\,f^*
P^{uvi_1\cdots i_{k-2}} \nonumber
\end{gather}
Using the equation above and the equations (\ref{comp}) we can
easily seen that the only possibility so that (\ref{assile}) holds
true is when and only when $\delta=\frac{2k-2+n}{1+n}.$

To prove (ii), denote by ${\mathfrak V}^\n$ the operators
(\ref{MultiSchwar2}) or (\ref{MultiSchwar2k}) written by means of
the connection $\n.$ Let $\tilde \n$ be another connection that is
projectively equivalent to $\n$ (see section \ref{compo}). We need
some ingredients for the proof. We will write the tensors
$\widetilde{\n}_u P^{i_1\cdots i_k}, \widetilde{\mathfrak{L}
(f)}_{ij}^v, \widetilde{\n}_u \,\widetilde{\mathfrak{L}
(f)}_{ij}^v,\n_i \,\widetilde{\mathfrak{L} (f)}_{j}$ and
$\widetilde{R_{ij}}$ in terms of $\n_u P^{i_1\cdots i_k},$ $
\mathfrak{L} (f)_{ij}^v,$ $ \n_u \,\mathfrak{L} (f)_{ij}^v,$ $\n_i
\,\mathfrak{L} (f)_{j},$ and $R_{ij}$ respectively. By using
(\ref{assoc}), we get
\begin{eqnarray}
\nonumber \widetilde{\n}_u P^{i_1\cdots i_k}&=&\n_u P^{i_1\cdots
i_k }+(2k-\delta(n+1))\, P^{i_1\cdots i_k}\,
\omega_u+\sum_{s=1}^k\,\delta^{i_s}_{u}\,\omega_v\,P^{vi_1\cdots
\widehat{i}_s\cdots i_k },
\end{eqnarray}
for all $P\in {\cal S}_\d^k(M)$, and
\begin{eqnarray}
 \nonumber \widetilde{\mathfrak{L}
(f)}_{ij}^v&=&\mathfrak{L} (f)_{ij}^v
+\mathrm{Sym}_{i,j}\,\delta^v_i\, {f^{-1}}^*\omega_j -
\mathrm{Sym}_{i,j}\,\delta^v_i\, \omega_j ,\\
\nonumber \widetilde{\n}_u \,\widetilde{\mathfrak{L}
(f)}_{ij}^v&=&\n_u \,\mathfrak{L}
(f)_{ij}^v+\mathrm{Sym}_{i,j}\,\delta^v_i\,
\nabla_u\,{f^{-1}}^*\omega_j-\mathrm{Sym}_{i,j}\,\delta^v_i\,
\nabla_u\,\omega_j-\mathrm{Sym}_{i,j}\,\omega_i\,
\widetilde{\mathfrak{L}(f)}_{uj}^v\\
&&+\delta^v_u\,\omega_t\, \widetilde{\mathfrak{L}(f)}^t_{ij}
-\omega_u\, \widetilde{\mathfrak{L}(f)}^v_{ij},
\end{eqnarray}
and finally
$$
\widetilde{R}_{ij}=R_{ij}+(n-1)\left ( \nabla_i\,
\omega_j-\omega_i\,\omega_j\right )\cdot
$$
By substituting these formul{\ae} into (\ref{MultiSchwar2}) we
obtain, after a long computation, that ${\mathfrak
V}^{\n}(f)={\mathfrak V}^{\tilde\n}(f).$

Suppose now $M=\bbR^n$ (or $M=\mathbb{S}^n$) and $M$ is endowed
with a projective structure. Let $f$ be a diffeomorphism belonging
to $\PSL(n+1,\mathbb{R}).$ Then there exist some constants $a_j^i,
b^i, c_l, d,$ where $i,j,l=1\ldots,n,$ such that
$$
f(x)=\left (\frac{a_j^1 x^j+b^1}{c_lx^l+d},\cdots, \frac{a_j^n
x^j+b^n}{c_lx^l+d}\right ).
$$ As the operators (\ref{MultiSchwar2}) and (\ref{MultiSchwar2k}) are
projectively invariant, we can take $\Gamma\equiv0.$ Therefore,
the tensor $\mathfrak T(f)^v_{ij}$ will take the form
$$
\mathfrak T(f)^v_{ij}=\frac{\partial^2 f^r}{\partial x^i\partial
x^j}\frac{\partial x^v}{\partial
f^r}-\frac{1}{n+1}\mathrm{Sym}_{i,j}\delta_i^v\, \frac{\partial^2
f^r}{\partial x^j\partial x^l}\frac{\partial x^l}{\partial f^r}
$$
It is a matter of a direct computation to prove that $\mathfrak
T(f)^v_{ij}\equiv 0,$ for all $f\in \PSL(n+~1,\mathbb{R}).$ Now,
directly from the equation (\ref{iya}) we see that ${\mathfrak
V}(f)\equiv0$ when $k=2.$ For $k>2,$ we will use again the
equation (\ref{ric}) and the proof is a long but straightforward
computation.\\
\cqfd
\subsection{A remark on the projective analogue of the
Laplace-Beltrami operator}
As a by-product of the formula (\ref{MultiSchwar2}) is the
projective analogue of the well-known Laplace-Beltrami operator
(see \cite{be}). It has been shown in \cite{bp1} that, for $k=2$
and for a particular value of $\delta,$ the conformal Schwarzian
derivative is given by the coboundary
$$
{f^*}^{-1} \Delta-\Delta,
$$
where $\Delta$ is the Laplace-Beltrami operator. In Theorem
(\ref{thm2}), we have proved that, for $\delta=\frac{n+2}{n+1},$
the projective Schwarzian derivative is the coboundary
$$
{f^*}^{-1}  B- B,
$$
where
$$
B:=\nabla_i\nabla_j-\frac{1}{n-1}\,R_{ij}.
$$
The operator $B$ is indeed projectively invariant; in virtue of
the conformal case, it can be then interpreted as the {\it
projective} analogue of the Laplace-Beltrami operator.
\subsection{Infinitesimal projective Schwarzian derivatives}
\begin{defi}
{\rm  (i) The infinitesimal operator associated with the operator
(\ref{Scw2}) is the operator
\begin{equation}
\label{liemul1}
 {\mathfrak
t}(X)\,(P)^{i_1\cdots i_s}:=\sum_{s=1}^{k-1}\left ({\mathfrak
l}(X)^{i_s}_{ij}-\frac{1}{n+1}\mathrm{Sym}_{i,j}\,
\delta_i^{i_s}\,{\mathfrak l}(X)_j\right ) \,P^{iji_1\cdots
\widehat{i}_s\cdots i_{k-1}},
\end{equation}
where ${\mathfrak l}$ is the 1-cocycle (\ref{mou}).

(ii) The infinitesimal operator associated with the operator
(\ref{MultiSchwar2}) and (\ref{MultiSchwar2k}) are respectively
the operators
\begin{gather} {\mathfrak u}(X)\,(P):= \displaystyle
\left ({\mathfrak l}(X)^{t}_{ij}-\frac{1}{n+1}\mathrm{Sym}_{i,j}\,
\delta_i^{t}\,{\mathfrak l}(X)_j\right )\,
\,\n_t\,P^{ij} \displaystyle +\n_t \,{\mathfrak l}(X)^{t}_{ij}\, P^{ij}\\
\phantom{{\mathfrak l}}{} -\frac{3+n-\delta(1+n)}{1+n}\nabla_i\,
{\mathfrak l}(X)_{j}\,P^{ij} +\frac{(1+n)(1-\delta)}{1+n}\,(L_X
R_{ij}) \,P^{ij},\label{liemul2p}
\end{gather}
and
\begin{gather}
{\mathfrak u}(X)\,(P)^{i_1\cdots i_{k-2}}:= \displaystyle
\sum_{s=1}^{k-2}\, \left ({\mathfrak
l}(X)^{i_s}_{ij}-\frac{1}{n+1}\mathrm{Sym}_{i,j}\,
\delta_i^{i_s}\,{\mathfrak l}(X)_j\right )\,
\,\n_t\,P^{ijti_1\cdots\widehat{i}_s\cdots i_{k-2} }\nonumber\\
\phantom{{\mathfrak u}(X)\,(P)}{}\displaystyle +\alpha_1\,\,\left
({\mathfrak
l}(X)^{t}_{ij}-\frac{1}{n+1}\mathrm{Sym}_{i,j}\,\delta_i^{t}\,{\mathfrak
l}(X)_j\right )\, \,\n_t\, P^{iji_1\cdots i_{k-2} }\label{liemul2}\\
\phantom{{\mathfrak u}(X)\,(P)}{} \displaystyle +\alpha_5\,\n_t
\,{\mathfrak l}(X)^{t}_{ij}\, P^{iji_1\cdots i_{k-2}}
+\alpha_2\,\sum_{s=1}^{k-2}\,\n_t\, {\mathfrak
l}(X)^{i_s}_{ij}\,P^{ijti_1\cdots\widehat{i}_s\ldots i_{k-2} }
\nonumber\\
\phantom{\!\!\!\!\!\!\!\!\!\!\!\!\!\!\!\!\!\!\!\!\!\!\!\!}{}+\alpha_6
\left (\nabla_i\, {\mathfrak l}(X)_{j}+\alpha_9\,L_XR_{ij}\right
)P^{iji_1\cdots i_{k-2}},\nonumber
\end{gather}
where the constants $\alpha_1,\alpha_2,\alpha_5,\alpha_6$ and
$\alpha_9$ are given as in (\ref{const}). }
\end{defi}
The following Corollaries result from Theorems (\ref{thm1}),
(\ref{thm2}) and (\ref{thm3}).
\begin{cor}
\label{coro1}
 For all $\delta\not=\frac{2k-1+n}{1+n},$ the map
$X\mapsto {\mathfrak t}(X)$ defines a non-trivial 1-cocycle valued
into ${\cal D}({\cal S}^{k}_\delta(M),{\cal S}^{k-1}_\delta(M)).$
Moreover, the operator (\ref{liemul1}) is projectively invariant,
namely it depends only on the projective class of the connection.
\end{cor}
\begin{cor}
\label{coro2} For all $\delta\not=\frac{2k-2+n}{1+n},$ the map
$X\mapsto {\mathfrak u}(X)$ defines a non-trivial 1-cocycle valued
into ${\cal D}({\cal S}^{k}_\delta(M),{\cal S}^{k-2}_\delta(M)).$
Moreover, the operators (\ref{liemul2p}) and (\ref{liemul2}) are
projectively invariant.
\end{cor}
\section{Conformally invariant Schwarzian derivatives}
Let $(M,\og)$ be a pseudo Riemannian manifold and let $\Gamma$ be
the Levi-Civita connection associated with the metric $\og.$
\begin{defi}
{\rm For all $f\in \Diff(M)$ and for all $P\in {\cal
S}^k_\delta(M),$ we put
\begin{equation}
\label{conf1} {\mathfrak A}(f)(P)^{i_1\cdots
i_{k-1}}=\mathrm{Coboundary}+c\sum_{s=1}^{k-1}\left(
\mathfrak{L}(f)_{ij}^{i_s}-\frac{1}{n}\,\mathrm{Sym}_{i,j}\,
\delta_i^{i_s}\,\mathfrak{L}(f)_j
\right)\,P^{iji_1\cdots\widehat{i}_s\cdots i_{k-1}},
\end{equation}
where $\mathfrak{L}(f)$ is the tensor (\ref{tenso}) and the
constant
$$
c=2-\delta n.
$$ }
\end{defi}
\begin{thm}\label{c1}
(i) For almost all values of $\delta$, the map $f\mapsto
{\mathfrak A} (f^{-1})$ defines a non-trivial 1-cocycle on
$\Diff(M)$ with values into $\End_{\mathrm{diff}} ({\cal
S}^k_{\d}(M),{\cal S}^{k-1}_\d(M))$;

(ii) The operator (\ref{conf1}) depends only on the conformal
class $[\og]$ of the metric. When $M=\bbR^n$  and $M$ is endowed
with a flat conformal structure, this operator vanishes on the
conformal group $\Og(p+1,q+1),$ where $p+q=n.$
\end{thm}
Now, we will introduce an other conformally invariant 1-cocycle
that takes values into $\End_{\mathrm{diff}}({\cal
S}_\delta^k(M),{\cal S}_\delta^{k-2}(M)).$ We suppose that $k>2;$
for $k=2,$ the 1-cocycle have already been introduced in
\cite{bp1}.

We denote by $R_{ij}$ the Ricci tensor components and by $R$ the
scalar curvature associated with the metric $\og.$
\begin{defi}
{\rm For all $f\in \Diff(M)$ and for all $P\in {\cal
S}^k_\delta(M),$ we put
\begin{gather}
\nonumber {\mathfrak B}(f)\,(P)^{i_1\cdots i_{k-2}}= \displaystyle
\mathrm{Coboundary}+\sum_{s=1}^{k-2}\,\left
(\mathfrak{L}(f)^{i_s}_{ij}-\frac{1}{n}\,\mathrm{Sym}_{i,j}\,
\delta^{i_s}_i\,\mathfrak{L}(f)_j \right)\n_t\,
P^{ijti_1\cdots\widehat{i}_s\cdots i_{k-2} }\\
\nonumber \phantom{\!\!\!\!\!\!
\!\!\!\!\!\!\!\!\!\!\!\!\!\!\!\!\!\!\!\!\!\!\!\!\!\!\!\!\!
\!\!\!\!\!\!\!\!}{} \displaystyle +\beta_1\,\left
(\mathfrak{L}(f)^{t}_{ij}-\frac{1}{n}\mathrm{Sym}_{i,j}\,\delta^{t}_i\,
\mathfrak{L}(f)_j \right) \,\n_t\,
P^{iji_1\cdots i_{k-2} }\\
\nonumber
\phantom{{\mathfrak B}(f)\,(P)}{}
+\left (\beta_2\,\n_t\,
\mathfrak{L}(f)^t_{ij}+\beta_3
\,\nabla_i\mathfrak{L}(f)_j+\beta_4\,\mathfrak{L}(f)_{i}\,
\mathfrak{L}(f)_{j}+
\beta_5\,\mathfrak{L}(f)_{ij}^u\,\mathfrak{L}(f)_{u}\right
)P^{iji_1\cdots i_{k-2} }\\
\phantom{{\mathfrak B}(f)}{} \displaystyle
+\sum_{s=1}^{k-2}\,\left (\beta_6\,\n_t\,
\mathfrak{L}(f)_{ij}^{i_s}+\beta_7\,\mathfrak{L}(f)^{i_s}_{ij}\,
\mathfrak{L}(f)_t+\beta_8\,\mathfrak{L}(f)_{ui}^{i_s}\,
\mathfrak{L}(f)_{jt}^u\,\right )
P^{ijti_1\cdots \widehat{i}_s\ldots i_{k-2}}\nonumber\\
\phantom{\!\!\!\!\!\!\!\!\!\!\!\!\!\!\!\!}{}
 \displaystyle +\left
(\beta_9\,({f^*}^{-1}R_{ij}-R_{ij})+\beta_{10}
({f^*}^{-1}R\,\og_{ij}-R\,\og_{ij})\right)\,
P^{iji_1\cdots i_{k-2}},\nonumber\\
\phantom{\!\!\!\!\!\!\!\!\!\!\!\!\!\!\!\!\!\!\!\!\!\!\!\!\!\! \!\!
\!\!\!\!\!\!\!\!\!\!\!\!\!\!\!\!\!\!\!\!\!\!\!\!\!\!}{}+e\displaystyle
\!\!\sum_{\substack{1\le s<t\le k-2}
}^{k-2}\,\mathfrak{L}(f)_{ij}^{i_s}\,\mathfrak{L}(f)_{pq}^{i_t}\,
P^{ijpqi_1\cdots \widehat{i}_s\cdots \widehat{i}_t\cdots i_{k-2}},
\label{conf2}
\end{gather}
where $\mathfrak{L}(f)_{ij}^k$ are the components of the tensor
(\ref{tenso}). The coefficient
$ e=\begin{cases} 1& \text{if $k\geq 4$},\\
0& \text{otherwise}
\end{cases},$ and the coefficients $\beta_1,\ldots,\beta_{10}$ are given
by
\begin{flalign*}
\beta_1&=\displaystyle \frac{1}{2}\,(4-2k+ n(\delta-1));& \beta_5
&=\displaystyle
\frac{1}{2}\,(1-\delta)(4-2k+n\,(\delta-1)); \\[3mm]
\beta_2 &= \displaystyle \frac{1}{2}\,(4-2k+ n(\delta-1))   ; &
\displaystyle \beta_6&=\displaystyle \frac{1}{6}\,(n+2k -\delta\,n);\\[3mm]
\beta_3 &
=\displaystyle \frac{1}{2}\,(\delta -1)(2-2k  +n(\delta-1));& \beta_7&
=(1-\delta);\\[3mm]
\beta_4&=\displaystyle \frac{1}{2}\,(\delta-1)^2 &
\beta_8&=\displaystyle \frac{1}{3}\,(6-2k+ n(\delta-1));\\
\end{flalign*}
$$
\begin{array}{cl}
\beta_9= &\displaystyle \frac{1}{6}\cdot\frac{4\,(6 - 5\, k + k^2)
-8\, n(k-2) (\delta-1)+3 n^2(\delta-1)^2}{n-2};\\[3mm]
\beta_{10}=&\displaystyle \frac{2(k-2) (2k(2k-5) + n(1 + 11\delta
-k(12\delta-7)) )-6n^3(\delta-1)^2 \delta}{12 (n-2) (n-1) (2 - 2 k
+ n (-1 +
2\delta))}\\[4mm]
&\displaystyle -\frac{ n^2(\delta-1)(2 +32\delta -
k(22\delta-5))}{12 (n-2) (n-1) (2 - 2 k + n (-1 + 2 \delta))}
\end{array}
$$
}
\end{defi}
\begin{thm}\label{c2}
(i) For almost values of $\delta$, the map $f\mapsto {\mathfrak B}
(f^{-1})$ defines a non-trivial 1-cocycle on $\Diff(M)$ with
values into $\cD ({\cal S}^k_{\d}(M),{\cal S}^{k-2}_\d(M))$;

(ii) The operator (\ref{conf2}) depends only on the conformal
class $[\og]$ of the metric. When $M=\bbR^n$ and $M$ is endowed
with a flat conformal structure, this operator vanishes on the
conformal group $\Og(p+1,q+1).$
\end{thm}
\subsection{The Algorithm and the proof of Theorems (\ref{c1}) and (\ref{c2})}
\label{fin3}
The operator
\begin{equation} \label {jah}
\sum_{s=1}^{k-1}\left
(\mathfrak{L}(f)^{i_s}_{ij}-\frac{1}{n}\,\mathrm{Sym}_{i,j}\,
\delta_i^{i_s}\mathfrak{L}(f)_j\right )
\end{equation}
satisfies obviously the 1-cocycle property. However, it lacks the
invariance property, in contradistinction with the operator
(\ref{ell}) which is projectively invariant. We will establish
here an Algorithm to transform the operator above into a
conformally invariant one.

Let us denote by $\cal C$ the 1-cocycle above written by means of
a connection $\nabla$ associated with the metric $\og$ and denote
by $\widetilde{\cal C}$ the same 1-cocycle written by means of a
connection belong to the same conformal class as described in
Section \ref{compo}. Using the formula (\ref{tof}), we get
$$
\left ( \widetilde{\cal C}(f)^{i_s}_{ij}-{\cal
C}(f)^{i_s}_{ij}\right ) P^{iji_1\cdots \widehat{i}_s\cdots
i_{k-1}} =\left
(\og_{ij}\,\,F^{i_s}-{f^*}^{-1}\og_{ij}\,{f^*}^{-1}\,F^{i_s}\right
)\,P^{iji_1\cdots \widehat{i}_s\cdots i_{k-1}}.
$$
In order to get ride the component
$\og_{ij}\,\,F^{i_s}\,P^{iji_1\cdots \widehat{i}_s\cdots
i_{k-1}}-{f^*}^{-1}\og_{ij}\,{f^*}^{-1}\,F^{i_s}\,P^{iji_1\cdots
\widehat{i}_s\cdots i_{k-1}},$  we adjust the 1-cocycle
(\ref{jah}) by incorporating the coboundary
$$
\gamma_1 ({f^*}^{-1}\, B_1-B_1),
$$
where $\gamma_1$ is a constant -- to be determined -- and
$B_1(P):=\sum_{s=1}^{k-1}\,\og_{uv}\,\og^{ti_s}\,\nabla_t\,
P^{uvi_1\ldots\widehat{i}_s \ldots i_{k-1}}.$

A direct computation using (\ref{tof}) proves that
$$
\begin{array}{cl}
\og_{uv}\,\og^{ti_s}\,\widetilde{\nabla}_t\,
P^{uvi_1\cdots\widehat{i}_s \cdots
i_{k-1}}=&\!\!\!\og_{uv}\,\og^{ti_s}\,\nabla_t\,P^{uvi_1\cdots\widehat{i}_s
\cdots i_{k-1}}+(k-\delta n)\,F^{i_s}\,\og_{uv}\,P^{uvi_1\cdots
\widehat{i}_s\cdots
i_{k-1}}\\[3mm]
&\displaystyle \!\!\!+\sum_{t=1, t\not= s}^{k-1}\!\!\og_{uv}\!\!
\left ( \og^{i_si_t}F_w\, P^{wuvi_1\cdots\widehat{i}_s
\cdots\widehat{i}_t \cdots i_{k-1}} \displaystyle
-F^{i_t}P^{uvi_1\cdots i_s\cdots\widehat{i}_t \cdots
i_{k-1}}\right )
\end{array}
$$
(For $k=2$ the last two terms will not to be taken into account.)

If we collect the coefficient of the component
$F^{i_s}\,\og_{uv}\,P^{uvi_1\cdots\widehat{i}_s \cdots i_{k-1}}$
and the component
${f^*}^{-1}\,F^{i_s}\,{f^*}^{-1}\,\og_{uv}\,P^{uvi_1\cdots\widehat{i}_s
\cdots i_{k-1}}$ we will get the equation
$$
\gamma_1 (2-\delta n)=c\cdot
$$
If $k=2, $ this is the only equation we need. In that case, the
coefficient $c$ and $\gamma_1$ are as in Table $1.$

\begin{center}
\label{tab1}
\begin{tabular}{c|cc|c}
 &  $c$ \rule[-5mm]{0mm}{12mm}  & $\gamma_1$ &
 $\mathfrak A$\\[2mm]\hline
$\displaystyle \delta =\frac{2}{n}$& $0$\rule[-5mm]{0mm}{12mm}
& $1$ & trivial\\[2mm]
$\displaystyle \delta \not=\frac{2}{n}$& $2-\delta n$ \rule[-5mm]
{0mm}{12mm}  & $1$& not trivial \\[2mm]\hline
\end{tabular}

\vspace*{0.5cm} Table 1.
\end{center}

If $k>2,$ the 1-cocycle
$$
\gamma_1 ({f^*}^{-1}\, B-B)+c \,\, {\cal C}(f)
$$
is still not conformally invariant. We have to incorporate, then,
another coboundary
$$
\gamma_2({f^*}^{-1}\, B_2-B_2),
$$
where $\gamma_2$ is a constant and
$B_2(P):=\sum_{s=1}^{k-2}\sum_{t\not=s}\og^{i_s
i_t}\,\og_{ij}\nabla_u\, P^{uiji_1\ldots \widehat{i}_s\ldots
\widehat{i}_t\ldots i_{k-1}}.$

Now, a direct computation using (\ref{tof}), we get
$$
\begin{array}{cl}
\displaystyle\og^{i_s i_t}\,\og_{ij}\widetilde{\nabla}_u\,
P^{uiji_1\cdots \widehat{i}_s\cdots \widehat{i}_t\cdots
i_{k-1}}=&\og^{i_s i_t}\,\og_{ij}\nabla_u\, P^{uiji_1\cdots
\widehat{i}_s\cdots
\widehat{i}_t\cdots i_{k-1}}\\[3mm]
&\displaystyle+(2k+n-4-\delta n)F_w\,\og^{i_si_t}\,\og_{ij}\,
P^{ijwi_1\cdots \widehat{i}_s\cdots
\widehat{i}_t\cdots i_{k-1}}\\[3mm]
&\displaystyle -\sum_{\substack{1\le l\le k-1\\l\not=
s,t}}\og^{i_si_t}\,\og_{uv}\,\og_{ij}\,F^{i_l}\,P^{ijuvi_1\cdots
\widehat{i}_l\cdots \widehat{i}_s\cdots\widehat{i}_t\cdots
i_{k-1}}
\end{array}
$$
(The last term should not be taken into account if $k=3.$)

Now, we collect the coefficient of the component
$F_m\,\og^{i_si_t}\,\og_{ij}\,P^{ijmi_1\cdots \widehat{i}_s\cdots
\widehat{i}_t\cdots i_k}$ we get the equation
$$
2\gamma_1+(2k-4+n(1-\delta))\gamma_2=0\cdot
$$
For $k=3,$ the coefficients $c, \gamma_1,$ and $\gamma_2$ are
given as in Table $2.$
\begin{center}
\label{tab2}
\begin{tabular}{c|ccc|c}
 &  $c$ \rule[-5mm]{0mm}{12mm}  & $\gamma_1$ & $\gamma_2$&
 ${\mathfrak A}$\\[2mm]\hline
$\displaystyle \delta =\frac{2}{n}$& $0$\rule[-5mm]{0mm}{12mm}   & $1$
&$\displaystyle-\frac{2}{n}$& trivial\\[2mm]
$\displaystyle \delta =\frac{2+n}{n}$& $0$ \rule[-5mm]{0mm}{12mm}
& $0$ &$1$& trivial \\[2mm]
$\delta$ not like above& $2-\delta n$ \rule[-5mm]{0mm}{12mm}  & $1$ &
$\displaystyle -\frac{2}{2+n(1-\delta)}$& not trivial\\[2mm] \hline
\end{tabular}

\vspace*{0.5cm} Table $2.$
\end{center}

If $k>3,$ the 1-cocycle
$$
\gamma_2 ({f^*}^{-1}\, B_2-B_2)+\gamma_1 ({f^*}^{-1}\,
B-B)+c\,\,{\cal C}(f)
$$
is still not conformally invariant. We have to incorporate then
another coboundary
$$
\gamma_3 ({f^*}^{-1}\, B_3-B_3)
$$ where
$B_3:=\sum_s\sum_{t}\sum_{p}\og^{i_si_t}\og^{i_pl}
\og_{uv}\,\og_{ab}\,\nabla_l P^{abuv i_1\cdots
\widehat{i}_s\cdots\widehat{i}_t\cdots\widehat{i}_u\cdots
i_{k-1}}.$ Then we proceed as before to find the constant
$\gamma_3.$ We will continue the procedure of incorporating
coboundaries up to the last coboundary:
$$
\gamma_k ({f^*}^{-1}\, B_k-B_k),
$$
where $B_k$ is an operator defined as follows:
\begin{enumerate}
\item If $k$ is even, then $B_k:=\mathrm{Sym}_{i_1,\ldots,
i_{k-1}}\,\og^{i_1i_2}\cdots \og^{i_{k-3}i_{k-2}}\og^{ui_{k-1}}
\og_{j_1j_2}\cdots \og_{j_{k-1}j_k}\,\nabla_u P^{j_1\cdots
j_{k}}.$ \item If $k$ is odd, then $B_k:=\mathrm{Sym}_{i_1,\ldots,
i_{k-1}}\,\og^{i_1i_2}\cdots \og^{i_{k-2}i_{k-1}}
\og_{j_1j_2}\cdots \og_{j_{k-2}j_{k-1}}\,\nabla_u P^{uj_1\cdots
j_{k-1}}.$

\end{enumerate}
The resulting 1-cocycle should be conformally invariant.

To prove that the operator (\ref{conf2}) satisfies the 1-cocycle
property is a long but straightforward computation using the
equations (\ref{comp}). It will determine the coefficients
$e,\beta_1,\ldots,\beta_8$ uniquely. In order to study the
invariance property, we need some ingredients. Using the relation
(\ref{tof}), we can prove that the following relations hold
$$
\begin{array}{cl}
\widetilde{\nabla}_u\,\widetilde{\mathfrak{L}(f)}^k_{ij}=&
\nabla_u\,\widetilde{\mathfrak{L}(f)}^k_{ij}-\mathrm{Sym}_{i,j}\,
F_j\,\widetilde{\mathfrak{L}(f)}^k_{iu}-F_u
\widetilde{\mathfrak{L}(f)}^k_{ij}\\[2mm]
&
\delta_u^k\,F_m\,\widetilde{\mathfrak{L}(f)}^m_{ij}+
\mathrm{Sym}_{i,j}\,\og_{uj}\,
F^t\,\widetilde{\mathfrak{L}(f)}^k_{it}-\og_{mu}\,F^k\,
\widetilde{\mathfrak{L}(f)}^m_{ij}\\[2mm]
\widetilde{\nabla}_uP^{i_1\cdots i_{k}}=&\nabla_uP^{i_1\cdots
i_{k}}+(k-\delta n)\,F_u\,P^{i_1\cdots i_{k}}+\displaystyle
\sum_{s=1}^{k}\delta^{i_s}_u\,F_m\,P^{mi_1\cdots
\widehat{i}_s\cdots i_k}\\[2mm]
& \displaystyle -\sum_{s=1}^k\og_{mu}\,F^{i_s}\,P^{mi_1\cdots
\widehat{i}_s\cdots i_k}\cdot
\end{array}
$$
Moreover,
$$
\begin{array}{cl}
\widetilde{R_{ij}}=& R_{ij}-(n-2)\left (\nabla_i
F_j-F_i\,F_j\right )-\left (\nabla_uF_v+(6-n)\,F_u\,F_v\right
)\og^{uv}\,\og_{ij}\\[2mm]
\widetilde{R}=& e^{-2F}\left
(R-(2n-2)\nabla_uF_v\,\,\og^{uv}-(7n-2-n^2)\,F_u\,F_v\,\og^{uv}\right
),
\end{array}
$$
where the wide tilde on each tensor means that the tensor is
written by means of a metric belonging to the conformal class.

In order to get a conformally invariant operator, we are required
to add the coboundary
$$
\mu_1 ({f^*}^{-1}\, B-B),
$$ where $B:=\og^{uv}\,\og_{ij}\nabla_u\nabla_v\,
P^{iji_1\cdots i_{k-2}}.$ Now, we proceed as above, to find the
constant $\beta_9$ and $\beta_{10}$ as well as the constant
$\mu_1$. We continue this process until we get a conformally
invariant operator.

\section{Schwarzian derivatives and Cohomology}
Let us first recall the following classical result (see \cite{ca,
wil}). Consider the space of Sturm-Liouville operators
$$A:=-2\,\frac{d^2}{dx^2}+u(x):\cF_{-\frac{1}{2}}
(\bbR)\rightarrow \cF_{\frac{3}{2}}(\bbR),$$ where $u(x)\in
\cF_2(\bbR)$ is the potential.

For all diffeomorphism $f\in \Diff(\bbR),$ the operator $f^*A$ is
still a Sturm-Liouville operator with potential $u\circ
f^{-1}\cdot {(f^{-1})'}^2+ S(f^{-1}),$ where $S(f^{-1})$ is the
Schwarzian derivative (\ref{lam}).

According to the Neijenhuis-Richardson's theory of deformation,
the space of Sturm-Liouville operators viewed as a
$\Diff(\bbR)$-module (~also as a $\Vect(\bbR)$-module) is a
non-trivial deformation of the quadratic differentials
$\cF_2(\bbR),$ generated by the Schwarzian derivative (see
\cite{bo1}).  More generally, the space of differential operators
acting on densities of arbitrarily weights is a non-trivial
deformation of a direct sum of densities of appropriate weights
(see \cite{bo1}). It is well-known that the problem of deformation
is related to the cohomology group
\begin{equation}
\label{coho1} \mathrm H^1(\Vect(\bbR),\Sl(2,\bbR);  \cD ({\cal
F}_\delta (\bbR),{\cal F}_{\delta'} (\bbR)))
\end{equation}
It has been proved in \cite{bo1} that the infinitesimal Schwarzian
derivative as well as other 1-cocycles generate this cohomology
group.
\begin{rmk}{\rm
The analogue cocycle on $\Vect(\bbR)$ associated with the
Schwarzian derivative is the so-called Gelfand-Fuchs cocycle:
$X\frac{d}{dx}\mapsto X'''\,dx^2$ (see e.g. \cite{f, gf}). }
\end{rmk}
Following these lines of though, we believe that, in
higher-dimension, the infinitesimal projective Schwarzian
derivatives are classes belonging to the cohomology group
$$
\mathrm H^1(\Vect(\bbR^n),\Sl(n+1,\bbR);  \cD ({\cal S}^k_\delta
(\bbR^n),{\cal S}^j_\delta (\bbR^n)))\cdot
$$
In the next section we will compute this cohomology group,
generalizing the result of \cite{lo2} for $\delta=0.$
\subsection{The projectively equivariant cohomology}
Consider $\mathbb{R}^n$ with the standard
$\SL(n+1,\mathbb{R})$-action as described in Section \ref{compo}.
\begin{thm}\label{mainth}
If $n>2,$ we have
\begin{gather}
\label{cal} \mathrm H^1(\Vect(\bbR^n),\Sl(n+1,\bbR);  \cD ({\cal
S}^k_\delta (\bbR^n),{\cal S}^j_\delta (\bbR^n)))=\nonumber\\[3mm]
\phantom{\mathrm H^1(\Vect(\bbR^n),\Sl(n+1,\bbR);  \cD ({\cal
S}^k_\delta (\bbR^n)}{}\left\{
\begin{array}{ll}
\bbR,& \mbox{if} \quad k-j=1,j\not=0
\mbox{ and }\,\,\delta\not=\frac{2k-1+n}{1+n},\\
\bbR,& \mbox{if} \quad k-j=2\,\,\mbox{ and }\,
\delta\not=\frac{2k-2+n}{1+n},\\
0,&\hbox{otherwise.}
\end{array}
\right.
\end{gather}
\end{thm}
The 1-cocycles that span the cohomology group above are the
operators (\ref{liemul1}), (\ref{liemul2p}) and (\ref{liemul2}).

The following remark will play a central r\^{o}le in our proof; it
has already been used in the papers \cite{bo1, lo2}. Let
$\mathfrak g$ be a Lie algebra, $\mathfrak h\subset \mathfrak g$
be a subalgebra and $M$ be a $\mathfrak g$-module. Any 1-cocycle
$c:\mathfrak g\rightarrow M$ that vanishes on the Lie sub-algebra
$\mathfrak h$ is automatically $\mathfrak h$-invariant. Indeed,
the 1-cocycle property reads
$$
L_X\,c(Y,A)-L_Y\,c(X,A)=c([X,Y],A)
$$
for all $X, Y\in \mathfrak g$ and for all $A\in M.$ Then
$$
L_X\,c(Y,A)=c([X,Y],A),
$$
which is nothing but the $\mathfrak h$-invariance property.

The strategy to proof Theorem (\ref{cal}) is as follows. We will
classify all $\Sl(n+1,\mathbb{R})$-invariant bilinear operators
from $\Vect(\mathbb{R})\otimes {\cal S}_{\delta}^k(\mathbb{R}^n)$
to ${\cal S}_{\delta}^j(\mathbb{R}^n)$, then we will isolate among
them 1-cocycles.
\subsubsection{$\Sl(n+1,\mathbb{R})$-invariant bilinear operators}
To begin with, we recall a lemma that has been proved in
\cite{lo2} for $\delta=0$ but the proof works well for any
$\delta.$
\begin{lem}
\label{le2} Every bilinear map from $\Vect(\mathbb{R})\otimes
{\cal S}_{\delta}^k(\mathbb{R}^n)$ to ${\cal
S}_{\delta}^j(\mathbb{R}^n)$ that is invariant with respect to the
action of the affine Lie algebra $\gl(n,\mathbb{R})\ltimes
\mathbb{R}^n$ is differentiable; moreover, it is given by the
divergence operator.
\end{lem}
{\bf Proof.} See \cite{lo2}
\begin{rmk} {\rm In fact, any 1-cocycle on $\Vect(M)$, where
$M$ is an arbitrary manifold, with values into $\cD ({\cal
S}_\delta^k(M),{\cal S}_\delta^j(M))$ is differentiable (cf.
\cite{lo2}). }
\end{rmk}
\begin{pro}
\label{class} The space of $\Sl(n+1,\mathbb{R})$-equivariant
bilinear operators form $\Vect(\mathbb{R})\otimes {\cal
S}_{\delta}^k(\mathbb{R}^n)$ to ${\cal
S}_{\delta}^{k-p}(\mathbb{R}^n)$ is as follows:

(i) for $k>p\geq 2,$ it is 2-dimensional;

(ii) for $k=p$, it is 1-dimensional;

(iii) for $p=1, k>2,$ it is 1-dimensional;

(iv) for $k=p=1,$ there is no such operators.
\end{pro}
{\bf Proof.} According to Lemma (\ref{le2}), any such operator
should have the expression
\begin{equation}
\label{you}
\begin{array}{cccl}
c(X,P)^{i_1\cdots i_{k-p}}&=&\displaystyle
\sum_{s=1}^p&\displaystyle \left ( \,\sum_{t=1}^{k-p}\alpha_s\,
\partial_{j_1}\cdots\partial_{j_{s+1}}(X^{i_t})\,
\partial_{j_{s+2}}\cdots \partial_{j_{p+1}}(P^{i_1\cdots
\widehat{i}_t\cdots i_{k-p}j_1\cdots j_{p+1}})\right.\\[3mm]
&&&\quad+\displaystyle \beta_s\,
\partial_{j_1}\cdots\partial_{j_{s}}\partial_k(X^{k})\,
\partial_{j_{s+1}}\cdots \partial_{j_{p}}(P^{i_1\cdots
i_{k-p}j_1\cdots j_{p}})\\[3mm]
&&&\quad+\displaystyle \left. \gamma_s\,
\partial_{j_1}\cdots\partial_{j_{s}}(X^{k})\,
\partial_k \partial_{j_{s+1}}\cdots \partial_{j_{p}}(P^{i_1\cdots
i_{k-p}j_1\ldots j_{p}})\right ),\\[3mm]
\end{array}
\end{equation}
where $\alpha_s,\beta_s,\gamma_s,$ for $s=1,\ldots,p,$ are real
numbers.

We will use the expression above as an Ansatz in order to classify
all $\Sl(n+1,\mathbb{R})$-invariant bilinear operators.

If we demand that the operator $c$ vanishes on the Lie algebra
$\Sl(n+1,\mathbb{R})$ we will impose the conditions
$$
\gamma_1=0,
$$
and
\begin{equation}
\label{mal}
 2\alpha_1 (k-p)+\beta_1 (1+n)+2\gamma_2=0\cdot
\end{equation}

A straightforward computation but quite complicated, prove that
the equivariance of the operator (\ref{you}) with respect to the
Lie algebra $\Sl(n+1,\mathbb{R})$ is equivalent to the following
system
\begin{eqnarray}
-s(s+2)\alpha_{s+1}+\gamma_{s+1}+(p-s)
(2k+n-p+s-\delta(1+n))\alpha_s&\!\!\!\!=&\!\!\!\!0,\label{e1}\\
-s(s+1)\beta_{s+1}+(s+1)\gamma_{s+1}
+(p-s)(2k+n-p+s-\delta(1+n))\beta_s&\!\!\!\!=&\!\!\!\!0,\\
-(s^2-1)\gamma_{s+1}+
(p-s)(2k+n-p+s-\delta(1+n))\gamma_s&\!\!\!\!=&\!\!\!\!0,\\
\label{eq1}
(s+1)\gamma_{s+1}+(k-p)(s+1)\alpha_{s}+(k-p+s-\delta(1+n))
\gamma_s+(1+n)\beta_s&\!\!\!\!=&\!\!\!\!0,
\end{eqnarray}
where $s=1,\ldots,p-1.$

The outcome (\ref{e1}) should not be taken into account if $k=p.$
\begin{lem}
For all $\delta,$ the system above is compatible.
\end{lem} {\bf Proof.} For $s=1$ the equation (\ref{eq1}) is nothing
but the equation (\ref{mal}). The proof follows by induction.\\
\cqf

Now we are ready to prove Proposition (\ref{class}).

(i) For $k>p\geq2,$ the space of solution is 2-dimensional spanned
by $\alpha_1, \beta_1.$

(ii) For $p=k,$ the constants $\alpha_s$ should be absence from
the system (\ref{e1}). The space of solution is 1-dimensional.

(iii) For $p=1,$ and $k>1,$ all the constants $\gamma_s$ are zero.
The space of solution is 1-dimensional generated by $\beta_1.$

(iv) For $k=p=1,$ all the constants $\gamma_s$ are zero and
$\beta_2$ should be absence form the equation (\ref{mal}). There
is no such operators.

(v) For $p=0,$ the space of solution is one-dimensional. \\
\cqf
\subsubsection{Proof of Theorem (\ref{mainth})}
The 1-cocycle property  of the operator (\ref{you}) adds to the
system above three other conditions:
\begin{equation}
\label{fem}
\begin{array}{lcl}
2(p-1)\,\alpha_1-\gamma_{p-1}&=&p\,\alpha_{p-1}\\
(p-1)\,\beta_1+\delta\,\gamma_2&=&\beta_2\\
\beta_1+\delta\,\gamma_p&=&\beta_p
\end{array}
\end{equation}
By using Proposition (\ref{class}), we get

(i) for $k>p\geq 2,$ we distinguish two cases.
\begin{enumerate}
\item If $p=2,$ the system above together with the condition
(\ref{fem}) admits (uniquely) a solution, independently on
$\delta.$ The corresponding cocycle associated with this class is
given in (\ref{liemul2}). This 1-cocycle turns into a trivial
cocycle for $\delta=\frac{2k-2+n}{1+n},$ as a consequence of
Corollary (\ref{coro2}). \item If $p>2,$ the system above together
with the conditions (\ref{fem}) admits (uniquely) a solution if
and only if
$$
\delta=\frac{2k-p+n}{1+n}\cdot
$$
\end{enumerate}
(ii) for $k=p.$ one distinguishes two cases:
\begin{enumerate}
\item If $k=p=2,$ the system above together with the condition
(\ref{fem}) admits (uniquely) a solution, independently on
$\delta.$ The corresponding 1-cocycle associated with this class
is given in (\ref{liemul2p}). This 1-cocycle turns into a trivial
cocycle for $\delta=\frac{2+n}{1+n},$ as a consequence of
Corollary (\ref{coro2}). \item if $k=p>2,$ the system above
together with the conditions (\ref{fem}) admits (uniquely) a
solution if and only if
$$
\delta=\frac{k+n}{1+n}\cdot
$$
\end{enumerate}

(iii) For $p=1,$ and $k>1,$ the unique 1-cocycle is given as in
(\ref{liemul1}). This 1-cocycle turns into a trivial cocycle for
$\delta=\frac{2k-1+n}{1+n},$ as a consequence of Corollary
(\ref{coro1}).

To achieve the proof of Theorem (\ref{mainth}) we are required to
prove the following Lemma.
\begin{lem}
For $\delta=\frac{2k-p+n}{1+n},$ any
$\Sl(n+1,\mathbb{R})$-invariant 1-cocycle from ${\cal S}_\delta^k
(\mathbb{R}^n)$ to ${\cal S}_\delta^{k-p} (\mathbb{R}^n)$ is
necessarily trivial.
\end{lem}
{\bf Proof.} The 1-cocycle conditions (\ref{fem}) turn the space
of solution of the system above into a 1-dimensional space. We are
led, then, to prove that any trivial 1-cocycle is necessarily
$\Sl(n+1,\mathbb{R})$-invariant for the particular value of
$\delta.$ To do that, we consider the operator $B$ defined as
follows. For all $P\in {\cal S}_\delta^k,$ we put
\begin{equation}
\label{bir1} B(P)=\partial_{j_1}\cdots \partial_{j_p} \,
P^{j_1\cdots j_pi_1\cdots i_{k-p}}\cdot
\end{equation}
Consider now the trivial 1-cocycle
\begin{equation}
\label{bir2} L_X\circ B-B\circ L_X.
\end{equation}
The order of the operator (\ref{bir2}) is $(p-1)$, because the
order of the operator (\ref{bir1}) is $p.$ One can easily seen
that the coefficients at any order less than $p-2$ contain
expressions in which the component $X$ is differentiated at least
three times. Thus, it vanishes on the Lie algebra
$\Sl(n+1,\mathbb{R}).$ Moreover, it is a matter of direct
computation to prove that the principal symbol of the operator
(\ref{bir2}) vanishes
on $\Sl(n+1,\mathbb{R})$ if and only if $\delta=\frac{2k-p+n}{1+n}.$\\
\cqf

Theorem (\ref{mainth}) is proven.\\
\cqfd
\subsection{Cohomology of $\Vect(M)$}
We need to recall the following Theorem.
\begin{thm}\cite{l}
\label{leco}
\begin{equation} \label{Liecal} \mathrm
H^1(\Sl(n+1,\bbR);  \cD ({\cal S}^k_\delta (\bbR^n),{\cal
S}^j_\delta (\bbR^n)))= \left\{
\begin{array}{ll}
\bbR,& \mbox{if} \quad k-j=0\\
\bbR^2,& \mbox{if} \quad k-j>0,\,\,\mbox{\rm and}\,\,\delta=\frac{k+j+n}{1+n}\\
0,&\hbox{\rm otherwise}
\end{array}
\right.
\end{equation}
\end{thm}
The 1-cocycles that span this cohomology group were given in
\cite{l}. These explicit expressions are as follows:
\begin{eqnarray}
\tau_j(X)(P)&=&\partial_i X^i\,\,
\partial_{l_1}\cdots\partial_{l_{k-j}}\,
P^{l_1\cdots l_{k-j}i_1\cdots i_{j}},\\
\kappa_j(X)(P)&=&\partial_{l_1}\partial_i X^i\,\,
\partial_{l_2}\cdots\partial_{l_{k-j}}\,
P^{l_1\cdots l_{k-j}i_1\cdots i_{j}},
\end{eqnarray}
for all $P\in {\cal S}_\delta^k(M).$

For $k-j=0,$ the cohomology group above is spanned by $\tau_k$.
For $k-j>0$ and $\delta=\frac{k+j+n}{1+n},$ it is spanned by
$\kappa_j$ and $\tau_j.$
\begin{pro}
\label{eng} (i) The 1-cocycles $\kappa_j$ can be extended uniquely
as 1-cocycles on $\Vect(\mathbb{R}^n)$ only for $k-j=1,2.$

(ii) The 1-cocycles $\tau_j$ can be extended uniquely to
$\Vect(\mathbb{R}^n)$ for $k-j=0,$ and for $(k,j)=(1,0)$ and
$\delta=1.$
\end{pro}
{\bf Proof.} (i) The 1-cocycles $\kappa_j$ can be extended to
$\Vect(\mathbb{R})$ for $k-j=1,2$ and the proof is just theirs
explicit expressions given in  (\ref{ext2p}), (\ref{ext2}) and
(\ref{ext3}). Let us prove the uniqueness. Suppose that there are
two 1-cocycles, say $c_1$ and $c_2,$ that extend $\kappa_j.$ This
implies that $c_1-c_2$ is zero on $\Sl(n+1,\mathbb{R}).$ The
1-cocycle $c_1-c_2$ is then projectively invariant. By using
Theorem (\ref{mainth}), the 1-cocycle $c_1-c_2$ should be a
coboundary, as $\delta=\frac{k+j+n}{1+n}.$ Thus, $c_1\equiv c_2.$

Now, we will proof that for $k-j>2,$ these 1-cocycles cannot be
extended. Suppose without loosing generality that $k-j=3.$ Any
1-cocycles that extend the 1-cocycles $\kappa_j$ should retains a
form as in (\ref{you}) but we incorporate another term
$\widehat{\beta}_1\,
\partial_i\,\partial_t\, X^t.$ The fact that the 1-cocycles in
question should coincide with the 1-cocycle $\kappa_j,$ leads to
the two conditions:
\begin{equation}
\label{sing} \gamma_1=0,\quad
2\gamma_2+(1+n)\,\beta_1+2(k-3)\,\alpha_1=0.
\end{equation}
The 1-cocycle property imposes the following conditions:
\begin{align*}
3\alpha_2-4\alpha_1+\gamma_2&=0,& 4\alpha_3-2\alpha_1+\gamma_3&=0,\\
6\alpha_3-3\alpha_2+\gamma_2&=0,&\beta_2-\delta\,
\gamma_2-2\,(\beta_1+\widehat{\beta}_1)&=0,\\
\beta_3-\delta\,\gamma_3-\beta_1-\widehat{\beta}_1&=0,&
3\beta_3-2\,\beta_2-(\delta-1)\,\gamma_2&=0,\\
3\gamma_3-2\,\gamma_2&=0,& (
2-\delta)(\beta_1+\widehat{\beta}_1)+(\delta-1)\,\beta_2&=0,\\
(\delta-1)\,\gamma_2+\beta_1+\widehat{\beta}_1&=0,&(\delta-1)\,
\alpha_2+\beta_1+\widehat{\beta}_1&=0
\end{align*}
The system above together with the outcomes (\ref{sing}) admits a
solution if and only if $\widehat{\beta}_1=0$ and
$\delta=\frac{2k-3+n}{1+n}.$ This means that the extended
1-cocycle is a coboundary and, moreover, vanishes on the Lie
algebra $\Sl(n+1,\mathbb{R}),$ which is absurd. This implies that
the 1-cocycle $\kappa_j$ cannot be extended. Part (i) is proven.

(ii) The 1-cocycles $\tau_j$ can be extended to
$\Vect(\mathbb{R})$ for $k-j=0$ and for $(k,j)=(1,0)$ and
$\delta=1.$ The proof is just theirs explicit expressions given in
(\ref{ext1}) and (\ref{ext4}). For the uniqueness, we can easily
proceed as in Part (i).

Suppose that the 1-cocycles $\tau_j$ can be extended to
$\Vect(\mathbb{R}^n)$ for the value of $k-j$ different from those
described above. Such 1-cocycles should retains a form as in
(\ref{you}) but we incorporate another term $\beta_0\,\partial_t\,
X^t.$ The fact that these 1-cocycles should coincide with the
1-cocycle $\tau_j$ once restricted to $\Sl(n+1,\mathbb{R}),$ leads
to the two conditions (\ref{sing}). Now, if we collect the
coefficient of the term
$\partial_iY^i\,\partial_{j_1}\,\partial_tX^t\partial_{j_2}\cdots
\partial_{j_p}$ we will get
$$
p\,\frac{k+j-1}{1+n}.
$$
This last outcome does not vanish, except when $(k,j)=(1,0),$ and
therefore $\delta=1.$ Part (ii) is proven. \\
\cqf

Let $M$ be any arbitrary manifold of dimension $n.$
\begin{thm}
\label{tah} For all $n>1,$ we have
\begin{equation} \label{curved} \mathrm
H^1(\Vect(M);  \cD ({\cal S}^k_\delta (M),{\cal S}^j_\delta (M)))=
\left\{
\begin{array}{ll}
\bbR\oplus {\mathrm H}_{\mathrm{DR}}^1(M),& \mbox{if} \quad k-j=0\\
\bbR,& \mbox{if} \quad k-j=1,j\not=0\\
\bbR^2\oplus {\mathrm H}_{\mathrm{DR}}^1(M),& \mbox{if} \quad
(k,j)=(1,0)\,\,\mbox{\rm and}\,\,\delta=1\\
\bbR,& \mbox{if} \quad k-j=2\\
0,&\hbox{\rm otherwise}
\end{array}
\right.
\end{equation}
\end{thm}
{\bf Proof.} For the proof we proceed as follows. Firstly, we
exhibit the 1-cocycles that span this cohomology group; secondly,
we proof the theorem for $\mathbb{R}^n$ then we extend the result
to an arbitrarily manifold.

(i) For $k-j=0,$ the 1-cocycles are already known (see \cite{f}).
\begin{gather}
\displaystyle \label{ext1} {\mathfrak a}_{\xi,\zeta}(X)\,(P)=
\displaystyle \left (\xi\,\Div(X)+ \zeta\, \omega (X)\right ) \,P,
\end{gather}
where  $\omega$ is a 1-form, $\Div(X)$ is the divergence operator
associated to some orientation and $\xi,\zeta$ are real numbers.

(ii) For $k-j=2,$ and $\delta \not =\frac{2k-2+n}{1+n},$ the
1-cocycle is given by the infinitesimal projective Schwarzian
derivative (\ref{liemul2}).

(iii) For $k-j=2,$ and $\delta=\frac{2k-2+n}{1+n},$ we distinguish
two cases:
\begin{enumerate}
\item For $k=2,$ the 1-cocycle in question is
\begin{gather}
\label{ext2p} \mathfrak{c}(X)(P)=\mathfrak{l}(X)_{i} \n_j P^{ij}+
\n_i\, \mathfrak{l}(X)_{j}\,P^{ij}\cdot
\end{gather}
where  ${\mathfrak l}(X)_{i}$ are the components of the tensor
(\ref{mou}).
\item For $k>2,$ the 1-cocycle is
\begin{equation}
\label{ext2}
\begin{array}{ccl} {\mathfrak c}(X)\,(P)^{i_1\cdots i_{k-2}}&=
&\displaystyle {\mathfrak l}(X)_{i}\,\nabla_j\,P^{iji_1\cdots
i_{k-2} }\\[3mm]
&&\displaystyle +\gamma_1\,\,\left ({\mathfrak
l}(X)^{t}_{ij}-\mathrm{Sym}_{i,j}\,\frac{1}{n+1}
\delta_i^{t}\,{\mathfrak l}(X)_j\right )\, \,\n_t\,
P^{iji_1\cdots i_{k-2} }\\[3mm]
&&\displaystyle +\gamma_2\,\n_t \,{\mathfrak l}(X)^{t}_{ij}\,
P^{iji_1\cdots i_{k-2}} +\gamma_3\,\sum_{s=1}^{k-2}\,\n_t\,
{\mathfrak l}(X)^{i_s}_{ij}\,P^{ijti_1\cdots\widehat{i}_s\cdots
i_{k-2} }
\\[3mm]
&&\displaystyle +\gamma_4 \nabla_i\, {\mathfrak
l}(X)_{j}\,P^{iji_1\cdots i_{k-2}},
\end{array}
\end{equation}
where the constants $\gamma_1,\ldots,\gamma_4$ are given by
\begin{equation}
\label{ferg}
\begin{array}{llll}
\nonumber \gamma_1&=\displaystyle \frac{1}{n+1};&\displaystyle
\gamma_2 &= \displaystyle \frac{1}{n+1}; \\[3mm]
\nonumber \gamma_3&=\displaystyle -\frac{1}{6}(1+n);&
\displaystyle \gamma_4 &=\displaystyle -\frac{1}{2}\,(2k-3).
\end{array}
\end{equation}
\end{enumerate}
(iv) For $k-j=1,$ $j\not=1$ and $\delta\not =\frac{2k-1+n}{1+n},$
the 1-cocycle is given by the infinitesimal projective Schwarzian
derivative (\ref{liemul1}).

(v)  For $k-j=1,$ $j\not=1$ and $\delta =\frac{2k-1+n}{1+n},$ the
1-cocycle is given by
\begin{gather}
\label{ext3} \displaystyle {\mathfrak c}(X)\,(P)^{i_1\cdots
i_{k-1}}= \displaystyle {\mathfrak l}(X)_{u}\, \,P^{ui_1\cdots
i_{k-1}},
\end{gather}
where  ${\mathfrak l}(X)_{i}$ are the components of the tensor
(\ref{mou}).

(vi) For $(k,j)=(1,0)$ and $\delta=1,$ the 1-cocycles are given by
\begin{gather}
\displaystyle \label{ext4} {\mathfrak d}_{\varepsilon,\xi,\zeta}
(X)\,(P)= \displaystyle \varepsilon \,{\mathfrak
l}(X)_i\,P^i+\left (\xi\,\Div(X)+ \zeta\, \omega (X)\right
)\nabla_i\, P^i.
\end{gather}
\subsubsection{Proof of Theorem (\ref{tah}) for the case
$M=\mathbb{R}^n$}
Let $c$ be any 1-cocycle on $\Vect(\mathbb{R}^n)$ with values into
$\cD ({\cal S}^k_\delta (\mathbb{R}^n),{\cal S}^j_\delta
(\mathbb{R}^n)).$ The restriction of this 1-cocycle, say
$\widehat{c},$ to $\Sl(n+1,\bbR)$ is obviously a 1-cocycle on
$\Sl(n+1,\bbR).$ We distinguish six cases:

(i) If $k-j>2,$ and $\delta\not=\frac{k+j+n}{1+n},$ then
$\widehat{c}$ is trivial, by Theorem (\ref{leco}). It follows that
there exists an operator, say $B,$ such that
$$
\widehat{c}(X)=[L_X,B],\quad \mbox{for all}\,\, X\in
\Sl(n+1,\bbR).
$$
Now, for all $X\in\Vect(\bbR^n)$ the map $X\mapsto c(X)-[L_X,B]$
is a 1-cocycle on $\Vect(\bbR^n)$ that vanishes on
$\Sl(n+1,\bbR).$ Theorem (\ref{mainth}) assures that such a
1-cocycle is trivial. A fortiori, $c\equiv0.$

(ii) If $k-j>2,$ and $\delta=\frac{k+j+n}{1+n},$ then
$\widehat{c}$ should be equal to zero by Proposition (\ref{eng}).
It implies that the 1-cocycle $c$ is vanishing on $\Sl(n+1,\bbR),$
and, thus, is trivial by Theorem (\ref{mainth}).

(iii) If $k=j,$ then $\widehat{c}$ is cohomologous to ${\mathfrak
a}_{1,0},$ by Theorem (\ref{leco}). It follows that there exists
an operator, say $B,$ such that
$$
\widehat{c}(X)-\alpha {\mathfrak c}_{1,0}(X)=[L_X,B], \quad
\mbox{for all}\,\, X\in \Sl(n+1,\bbR).
$$
Now, for all $X\in\Vect(\bbR^n)$ the map $X\mapsto c(X)-\alpha
{\mathfrak a}_{1,0}(X)-[L_X,B]$ is a 1-cocycle on $\Vect(\bbR^n)$
that vanishes on $\Sl(n+1,\bbR).$ Theorem (\ref{mainth}) assures
that such a 1-cocycle is necessarily trivial. A fortiori, $c\equiv
{\mathfrak a}_{1,0}.$

(iv) If $k-j=1,$ and $j\not=1,$ we will prove that $c$ is
cohomologous to one of the 1-cocycles (\ref{ell}) or (\ref{ext3}),
depending on the value of $\delta.$
\begin{enumerate}
\item For $\delta\not =\frac{2k-1+n}{1+n}$, the 1-cocycle $\widehat{c}$
should be trivial by Theorem (\ref{leco}). It follows that there
exists an operator, say $B,$ such that
$$
\widehat{c}(X)=[L_X,B], \quad \mbox{for all}\,\, X\in
\Sl(n+1,\bbR).
$$
Now, for all $X\in\Vect(\bbR^n)$ the map $X\mapsto c(X)-[L_X,B]$
is a 1-cocycle on $\Vect(\bbR^n)$ that vanishes on
$\Sl(n+1,\bbR).$ Theorem (\ref{mainth}) assures that such a
1-cocycle is necessarily unique. A fortiori, $c\equiv {\mathfrak
t}.$ \item For $\delta=\frac{2k-1+n}{1+n}$, the 1-cocycle
$\widehat{c}$ should be cohomologous to the 1-cocycle $\alpha
\kappa_{k-1}+\beta\tau_{k-1}$, by Theorem (\ref{leco}). Moreover,
by using proposition (\ref{eng}) the 1-cocycle $\kappa_{k-1}$ is
the only 1-cocycle that can be extended. It follows that there
exists an operator, say $B,$ such that
$$
\widehat{c}(X)-\alpha \kappa_{k-1}(X)=[L_X,B], \quad \mbox{for
all}\,\, X\in \Sl(n+1,\bbR).
$$
Now, for all $X\in\Vect(\bbR^n),$ the map $X\mapsto
c(X)-{\mathfrak c}(X)-[L_X,B]$ is a 1-cocycle on $\Vect(\bbR^n)$
that vanishes on $\Sl(n+1,\bbR).$ Theorem (\ref{mainth}) assures
that such a 1-cocycle is necessarily trivial. A fortiori, $c\equiv
{\mathfrak c}.$
\end{enumerate}

(v) If $(k,j)=(1,0)$ and $\delta=1.$ By using the same method as
before, we can prove that $c$ is cohomologous to the 1-cocycles
${\mathfrak d}_{\varepsilon,\xi,0}.$

(vi) If $k-j=2,$ By using the same method as before, we can prove
that $c$ is cohomologous to the 1-cocycles (\ref{liemul2p}) or
(\ref{ext2}).

Theorem (\ref{tah}) is proven for $\mathbb{R}^n.$
\subsubsection{Proof of Theorem (\ref{tah}) for the case of an arbitrary
manifold}
The techniques that we are going to use here have been already
used in \cite{lo2} for $\delta=0.$

(i) For $k-j=0$ we have
$$\mathrm{H}^1(\Vect(M);{\mathcal
D}({\mathcal S}^k_\delta(M),{\mathcal S}^k_\delta(M)))\simeq
\mathrm{H}^1(\Vect(M);C^{\infty}(M)).$$ The later cohomology group
is well-known; it is isomorphic to $\mathbb{R}\oplus {\mathrm
H}_{\mathrm{DR}}^1(M)$ (see, e.g., \cite{f}).

(ii) For $(k,j)=(1,0)$ and $\delta=1$ we have
$$
\mathrm{H}^1(\Vect(M);{\mathcal D}({\mathcal S}^1_1(M),{\mathcal
S}^0_1(M)))\simeq \mathrm{H}^1(\Vect(M);\Omega^1(M))\oplus
\mathrm{H}^1(\Vect(M);C^{\infty}(M)).$$ For the proof we proceed
as follows. Let $c$ be a 1-cocycle on $\Vect(M)$ with values into
${\mathcal D}({\mathcal S}^1_1(M),{\mathcal S}^0_1(M)).$ The fact
that $M$ is endowed with a connection implies that the 1-cocycle
$c$ can be written as
$$
b(X)\nabla_i \,+a_i(X)\cdot
$$
The 1-cocycle condition of the 1-cocycle $c$ implies that the
components $a_i$ should define a 1-cocycle belonging to the
cohomology group $\mathrm{H}^1(\Vect(M);\Omega^1(M))$ and $b$
should define a 1-cocycle belonging to the cohomology group
$\mathrm{H}^1(\Vect(M);C^{\infty}(M)).$ Reciprocally, any two
1-cocycles in $\mathrm{H}^1(\Vect(M);\Omega^1(M))$ and
$\mathrm{H}^1(\Vect(M);C^{\infty}(M))$ will define the 1-cocycle
$c,$ as it is given above. The cohomology group
$$
\mathrm{H}^1(\Vect(M);\Omega^1(M))
$$ is well-known; it is isomorphic to $\mathbb{R}$
(see, e.g., \cite{tsu}).

(iii) For $k-j>2.$ Let $c$ be a 1-cocycle on $\Vect(M)$ valued
into ${\mathcal D}({\mathcal S}^k_\delta(M),{\mathcal
S}^j_\delta(M)).$ On a local chart $U,$ the restriction $c_{|_U}$
is trivial. Namely, it exists an operator, say $B_{|_U},$ on $U$
such that
$$
c_{|_U}=L_X(B)_{|_U}\cdot
$$
A local coordinates patching will be used to extend the operator
$B_{|_U}.$ To do that, we should prove that $B_{|_U}=B_{|_V}$ on
the intersection $U\cap V.$ Indeed,
$$
0=c_{|_{U\cap V} }-c_{|_{U\cap
V}}=L_X(B)_{|_U}-L_X(B)_{|_V}.
$$
As there is no $\Vect(M)$-invariant operators for $k-j>2,$ it
implies that $B_{|_U}=B_{|_V}$ on $U\cap V.$

(iv) For $k-j=2$ and $k>2.$ Let $c$ be a 1-cocycle on $\Vect(M)$
valued into ${\mathcal D}({\mathcal S}^k_\delta(M),{\mathcal
S}^j_\delta(M)).$ On a local chart $U,$ the restriction $c_{|_U}$
is cohomologous to the 1-cocycle (\ref{liemul2p}) or (\ref{ext2}).
Namely, it exists an operator, say $B_{|_U},$ on $U$ such that
$$
c_{|_U}+\alpha_{U} {\mathfrak p}(X)=L_X(B)_{|_U}\cdot
$$
where $\mathfrak p$ is one of the two 1-cocycles (\ref{liemul2p})
or (\ref{ext2}). On the intersection $U\cap V,$ one has
$$
(\alpha_U-\alpha_V)\,{\mathfrak p}(X)=L_X(B)_{|_U}-L_X(B)_{|_V}.
$$
Thus, $\alpha_U-\alpha_V=0$ because ${\mathfrak p}$ is not a
coboundary and, a fortiori, $B_{|_U}=B_{|_V}$ on $U\cap V,$ as
there is no $\Vect(M)$-invariant operators for $k-j=2.$

(v) For $k=2$ and $j=0,$ the proof is the same as in (iii).

(vi) For $k-j=1,$ and $j\not=1$, the proof is the same as in
(iii).
\subsection{Cohomology of $\Diff(S^n)$}
In order to compute the cohomology of the group of diffeomorphisms
$\Diff(M),$ we deal with differential cohomology ``Van Est
Cohomology''; this means we consider only differential cochains
(see~\cite{f}). The more general case -- namely, the cohomology
with also non-differentiable cochains -- is an intricate problem,
and even though no explicit cocycles are known in our situation.

 Let $\mathbb{S}^{n}$ be the $n$-dimensional sphere. It is
 well-known that the maximal compact group of ``rotations''
of $\mathbb{S}^n$, $SO(n+1),$ is a deformation retract of the
group $\mathrm{Diff}_+ (\mathbb{S}^n),$ for $n=1,2,3$
(see~\cite{th}). Since the space $\mathrm{Diff}_+ (S^n)/SO(n+1)$
is acyclic, the Van Est cohomology of the Lie group
$\mathrm{Diff}_+(\mathbb{S}^n)$ can be computed using the
isomorphism (see, e.g., \cite[p.~298]{f})
\begin{equation}
\label{bah} {\mathrm H}^1(\mathrm{Diff}_{+}(\mathbb{S}^n);
{\mathcal D}({\mathcal S}^k_\delta (\mathbb{S}^n),{\mathcal
S}^j_\delta (\mathbb{S}^n)))\simeq \mathrm{H}^1(\mathrm{Vect}
(\mathbb{S}^n), SO(n+1);{\mathcal D}({\mathcal S}^k_\delta
(\mathbb{S}^n),{\mathcal \mathbb{S}}^j_\delta(\mathbb{S}^n))).
\end{equation}
We state the following Theorem that generalizes the result of
\cite{bn} for $\delta=0.$
\begin{thm}
\label{sala} For $n=2,3,$ the first-cohomology group
\begin{equation}
\label{mainlast} \mathrm{H}^1(\mathrm{Diff}_{+}(\mathbb{S}^n);
{\mathcal D}({\mathcal S}^k_\delta(\mathbb{S}^n),{\mathcal
S}^j_\delta(\mathbb{S}^n)))= \left\{
\begin{array}{ll}
\bbR,& \mbox{if} \quad k-j=0\\
\bbR,& \mbox{if} \quad k-j=1,j\not=0\\
\bbR^2,& \mbox{if} \quad (k,j)=(1,0)\,\,\mbox {\rm and}\,\,\delta=1\\
\bbR,& \mbox{if} \quad k-j=2\\
0,&\hbox{\rm otherwise}
\end{array}
\right.
\end{equation}
\end{thm}
{\bf Proof.} We will first give the explicit 1-cocycles that span
the cohomology group above.

(i) For  $k-j=0$. Any diffeomorphism $f\in
\mathrm{Diff}_{+}(\mathbb{S}^n)$ preserves the volume form on
$\mathbb{S}^n$ up to some factor. The logarithm function of this
factor defines a 1-cocycle on $\mathrm{Diff}(\mathbb{S}^n),$ say
${\mathfrak J}(f),$ with values in $C^{\infty}(\mathbb{S}^n).$
Now, the 1-cocycle in question is just the multiplication operator
by ${\mathfrak J}(f).$

(ii) For  $k-j=1,$ $j\not=0$ and $\delta\not=\frac{2k-1+n}{1+n},$
the 1-cocycle in question is the Schwarzian derivative
(\ref{Scw2}). For $\delta=\frac{2k-1+n}{1+n},$ the 1-cocycle is
$$
\mathfrak{L}(f)_u\,P^{ui_1\cdots i_{k-1}}
$$
where $\mathfrak{L}(f)_u$ are the components of the trace of the
tensor (\ref{ell}).

(iii) For $(k,j)=(1,0)$ and $\delta=1$ the 1-cocycles are
\begin{gather}
\displaystyle \nonumber \varepsilon
\mathfrak{L}_u(f)\,P^{u}+\xi\,{\mathfrak J}(f)\nabla_u\,P^{u},
\end{gather}
where $\varepsilon$ and $\xi$ are real numbers.

 (iv) For  $k-j=2$ and $\delta\not=\frac{2k-2+n}{1+n}$, the
 1-cocycle is the Schwarzian derivative (\ref{MultiSchwar2k}).
For $\delta=\frac{2k-2+n}{1+n}$, we distinguish two cases:
\begin{enumerate}
\item For $k=2,$ the 1-cocycle in question is
\begin{gather}
\label{perm} \mathfrak{b}(f)(P)=\mathfrak{L}(f)_{i} \n_j P^{ij}+
\n_i\, \mathfrak{L}(f)_{j}\,P^{ij}-\frac{1}{2(n+1)}\,
\mathfrak{L}(f)_{i}\mathfrak{L}(f)_{j}P^{ij}\cdot
\end{gather}
\item For $k>2,$ the 1-cocycle is
\begin{gather} \nonumber
{\mathfrak b}(f)\,(P)^{i_1\cdots i_{k-2}}= {\mathfrak L}(f)_{u}\,
\,\n_{v}\,P^{tuvi_1\cdots i_{k-2} } +\beta_1\,\,{\mathfrak
T}(f)^t_{uv}\, \,\n_t\,P^{uvi_1\cdots i_{k-2} }\\
\quad {} \nonumber +\sum_{s=1}^{k-2}\left (\beta_2\n_{t}\,
\mathfrak{L}(f)_{uv}^{i_s}+\beta_3\,
\mathfrak{L}(f)_{wt}^{i_s}\,\mathfrak{L}(f)_{uv}^w \right )
P^{tuvi_1\cdots \widehat{i}_s\cdots i_{k-2} }\\
\quad {}\nonumber+\left
(\beta_4\,\mathfrak{L}(f)_{uv}^w\,\mathfrak{L}(f)_{w}+\beta_5\,\n_t
\mathfrak{L}(f)^t_{uv}+\beta_6 \,\nabla_u\, \mathfrak{L}(f)_v
\right ) \,P^{uvi_1\cdots i_{k-2}},
\end{gather}
where $\mathfrak{L}(f)^k_{ij}$ are the components of the tensor
(\ref{tenso}) and  ${\mathfrak T}(f)^k_{ij}$ are the components of
the tensor (\ref{ell}). The constants $\beta_1,\ldots,\beta_6$ are
given by
\begin{equation}
\label{ferg1}
\begin{array}{llllll}
\nonumber \beta_1&=\displaystyle \frac{1}{n+1};&
\beta_3&=\displaystyle \frac{1}{3}(1+n) ;&\displaystyle
\beta_5 &= \displaystyle \frac{1}{n+1}; \\[3mm]
\nonumber \beta_2&=\displaystyle -\frac{1}{6}(1+n); &
\beta_4&=\displaystyle \frac{1}{2}\,(2k-3);& \displaystyle \beta_6
&=\displaystyle -\frac{1}{2}\,(2k-3).
\end{array}
\end{equation}
\end{enumerate}
We are now ready to prove Theorem (\ref{sala}). First, observe
that the De Rham classes in the cohomology group (\ref{cal}) is
trivial since $\mathrm{H}^1_{\mathrm {DR}}(\mathbb{S}^n)=0.$

In view of (\ref{bah}), the cohomology group
$\mathrm{H}^1(\mathrm{Diff}_{+}(\mathbb{S}^n); {\mathcal
D}({\mathcal S}^k_\delta(\mathbb{S}^n),{\mathcal
S}^j_\delta(\mathbb{S}^n)))=0$ for $k-j\neq 0, 1, 2,$ and for
$(k,j)=(1,0)$ with $\delta\not =1.$ Besides,

(i) For $k-j=0,$ suppose that there are two 1-cocycles
representing cohomology classes in the cohomology
group~(\ref{mainlast}). The isomorphism above shows that these two
1-co\-cycles induce two non-cohomologous classes in the cohomology
group (\ref{cal}), which is absurd.

(ii) For $k-j=1,$ and $j\not=0$, idem.

(iii) For $k-j=2$, idem.

(iv) For $(k,j)=(1,0)$ and $\delta=1,$ suppose that there are more
than two 1-cocycles representing cohomology classes in the
cohomology group~(\ref{mainlast}). The isomorphism above shows
that these 1-co\-cycles induce non-cohomologous classes in the
cohomology group (\ref{cal}), which is absurd.

Theorem~\ref{sala} follows, therefore, from explicit constructions
of the 1-cocycles above.
\begin{rmk}
{\rm Theorem (\ref{sala}) remains true as far as the rotation
group $SO(n+1)$ is a deformation retract of the group
$\Diff_+(\mathbb{S}^n)$ for all $n.$ We do not know whether this
statement is true or not.

}
\end{rmk}
\subsection{Relation to the Vey Cocycle}
Throughout this section, we will assume that $\delta=0$. The main
result is to give a relation between the projective Schwarzian
derivative (\ref{MultiSchwar2k}) and the well-known Vey cocycle,
answering a question raised in \cite{bn}.

Recall that the Vey cocycle is a object that is closely related to
deformation quantization (see \cite{v} for more details.). It is,
in fact, a cohomology class that span the component $\mathbb{R}$
of the cohomology group ${\mathrm
H}^2(C^{\infty}(T^*M),C^{\infty}(T^*M))\equiv {\mathrm
H}^2_{\mathrm{DR}}(M)\oplus \mathbb{R}$ (see \cite{v}). In order
to write it down, we need to lift the connection to a connection
on the cotangent bundle $T^*M$ (see \cite{yano} for more details).
We are mainly interested when its first component is restricted to
$\Vect(M)\subset C^{\infty}(T^*M).$ The Vey cocycle reads
accordingly as follows.
\begin{gather}
S^3(X):=Sym_{\mathbf j,\mathbf i,\mathbf k} \left ({\mathfrak l
}(\tilde X)_{\mathbf m \mathbf l}^{\mathbf j} \!\cdot\!
\mathrm{\omega}^{\mathbf{im} }\! \cdot \mathrm{\omega}^{
\mathbf{kl} } \right ) \tilde \nabla_{\mathbf i}  \tilde
\nabla_{\mathbf j} \tilde \nabla_{\mathbf k}\label{exp32}.
\end{gather}
In the formula above, the subscript $\mathbf i$ runs from $1$ to
$2n,$ and $\mathrm{\omega}$ stands for the standard symplectic
structure on $T^*M,$ and $\tilde X$ is the Hamiltonian lift of
$X.$

The following cocycle were introduced in \cite{bn}, and
interpreted as a {\it group Vey cocycle}:
\begin{gather}
GS^3(f):=Sym_{\mathbf j,\mathbf i,\mathbf k} \left ({\mathfrak
L}_{\mathbf m \mathbf l}^{\mathbf j}(\tilde f) \!\cdot\!
\mathrm{\omega}^{\mathbf{im} }\! \cdot \mathrm{\omega}^{
\mathbf{kl} } \right ) \tilde \nabla_{\mathbf i}  \tilde
\nabla_{\mathbf j} \tilde \nabla_{\mathbf k}\nonumber\\
\phantom{{\mathcal L}(f):=}{}-\frac{3}{2} Sym_{ \mathbf n,\mathbf
m,\mathbf i} \left ( {\mathfrak L}_{\mathbf l \mathbf k}^{ \mathbf
n}(\tilde f)\! \cdot\! \mathrm{\omega}^{ \mathbf{ml} } \!\cdot\!
\mathrm{\omega}^{ \mathbf{ik} }\right ) \cdot {\mathfrak
L}_{\mathbf {mn}}^{\mathbf {j}}(\tilde f) \tilde \nabla_{\mathbf
i} \tilde \nabla_{\mathbf j},\label{exp33}
\end{gather}
where $\tilde f$ is the symplectic lift of $f$ to $T^*M$ and
${\mathfrak L}(f)^k_{ij}$ are the components of the tensor
(\ref{tenso}) with respect to the lifted connection on $T^*M.$
\begin{pro} \label{fina}
The relation between the Vey cocycle and the
projective Schwarzian derivative is as follows:

(i) For all $X\in\Vect(M),$ we have
$$
\begin{array}{ccl}
{\mathfrak v}(X)^{i_1\cdots i_{k-2}}&=&\displaystyle
\frac{1}{2}\,L_X (\nabla_i
\,\nabla_j) +\frac{2-2k-n}{2}\,S^3(X)_{|_{{\cal S}^k(M)}}\\[3mm]
&& \displaystyle
+\frac{11+4k^2-2n(5-4k)+3n^2-12k}{6-6n}\,\,L_X(R_{ij}).
\end{array}
$$
(ii) For all $f\in \Diff(M),$ we have
\begin{equation}
\label{fina2}
\begin{array}{ccl}
{\mathfrak V}(f)^{i_1\cdots i_{k-2}}&=&\displaystyle
\frac{1}{2}\,{f^{-1}}^*(\nabla_i\, \nabla_j)-\nabla_i \, \nabla_j
+\frac{2-2k-n}{2}\,GS^3(f)_{|_{{\cal
S}^k(M)}}\\[3mm]
&& \displaystyle
+\frac{11+4k^2-2n(5-4k)+3n^2-12k}{6-6n}\,\,({f^{*}}^{-1}R_{ij}-R_{ij}).
\end{array}
\end{equation}
\end{pro}
{\bf Proof.} For the proof, we have to expound the formulas
(\ref{exp33}) and (\ref{exp32}) once restricted to ${\cal S}^k(M)$
and write these expressions in terms of the initial connection on
$M.$ Then, the proof follows by a direct computation.
\subsection{Conclusion and Open Problems}
The programm for defining the projective and conformal
multi-dimensional Schwarzian derivatives is achieved now in this
paper. However, it would be interesting to investigate topological
properties of these derivatives. For instance, it has recently
been proved that the classical Schwarzian derivative of a
diffeomorphism admits at least four zeros in \cite{ot}. According
to Ghys-Ovsienko-Tabachnikov, this property is the {\it four
vertex Theorem} of a time-like curve on the torus endowed with a
Lorentzian metric. It would be interesting to know whether a
theorem of this type holds true for our multi-dimensional
Schwarzian derivatives.

According to Theorem (\ref{sala}), the conformal Schwarzian
derivatives are only the operators (\ref{conf1}) and
(\ref{conf2}), except another cocycle may appear for the
particular values $(k,j)=(1,0)$ and $\delta=1.$ But, we do not
expect new cocycles other than those given here. More precisely,
we are led to compute the cohomology group
$$
\mathrm{H}^1(\Diff(\mathbb{R}^n), \Og(p+1,q+1); {\mathcal
D}({\mathcal S}^k_\delta(\mathbb{R}^n),{\mathcal
S}^j_\delta(\mathbb{R}^n))\cdot
$$
The computation of this cohomology group is more intricate, and
even though for the cohomology of $\Vect(\mathbb{R}^n)$ the
computation is still out of rich.

The conformal Schwarzian derivative is certainly related to the
Vey cocycle and an analogue to the Proposition (\ref{fina}) is
certainly true. We are required to incorporate to the Vey cocycle
an appropriate coboundary to get a formula analogous to that in
(\ref{fina2}). We recall that this coboundary has been added, as
explained in section (\ref{fin3}), in order to get the invariance
property.

Recently, the author has investigated an analogue of the operator
(\ref{conf1}) to the ({\it generic}) Finsler structures in
\cite{fin}, using some connections associated with the Finsler
structure. This operator has the property that it coincides with
the operator (\ref{conf1}) when the Finsler structure is
Riemannian. It would be interesting to investigate Schwarzian
derivatives in other geometry; for instance: CR structures,
quaternionic structures...

It should be stressed that in the literature alternative
approaches were developed in order to extend the classical
Schwarzian derivative to a multi-dimensional manifold (see for
example
\cite{a, c, g, mm, os, ol, sa1, sa2}).\\
\\
\noindent {\it Acknowledgement.} The problem of investigating
invariant Schwarzian derivatives is an idea due to V. Ovsienko
proposed as a subject for my thesis and supervised by him. I am
grateful to him for his constant support.


\end{document}